\documentclass{amsart}
\usepackage[T1]{fontenc}
\usepackage{amssymb}
\usepackage{amscd}
\usepackage[hyperref]{degt}

\def\2{\color{red}}
\def\3{\color{magenta}}

\title{Real plane sextics without real points}

\def\aset#1{\xdef\lastset{#1}%
 \smash{\raise9pt\vbox{\hypertarget{#1}{}}}{\singset{#1}}}
\def\sset#1{\hyperlink{#1}{\singset{#1}}}
\def\same{\omit\quad---\hss}

\makeatletter
\def\get@d#1{\xdef\ppd{#1}}
\def\get@count#1{\xdef\ppcount{#1}}
\def\get@rcount#1{\xdef\pprcount{#1}}
\def\get@conics#1{\xdef\ppconics{#1}}
\def\get@forms#1{\xdef\ppforms{#1}}
\def\get@D#1{\xdef\ppD{#1}}
\def\get@A#1{\xdef\ppA{#1}}
\def\get@sym#1{\xdef\ppsym{#1}}
\def\get@all{%
 \let\d\get@d
 \let\count\get@count
 \let\rcount\get@rcount
 \let\conics\get@conics
 \let\forms\get@forms
 \let\D\get@D
 \let\A\get@A
 \let\sym\get@sym
 \A{}\D{}\global\let\sep\relax\global\let\comma;
}
\def\DA#1#2{\hyperref[DA-#2]{#1^{\mathrm{\romannumeral#2}}}}
\def\str#1{$\DA{#1}1$}
\def\sstr#1{$\DA{#1}2$}
\def\rsmash#1{\setbox\z@\hbox to\z@{$\m@th#1$\hss}\ht\z@\z@\box\z@}
\def\printDA{\setbox\z@\hbox{\let\comma\relax$\m@th
 \ifx\ppD\empty\else\DA{\bD_6}\ppD\let\comma,\fi
 \ifx\ppA\empty\else\comma\DA{\bA_3}\ppA\fi$}%
 \ifdim\wd\z@=\z@\else(\ht\z@\z@\box\z@)\let\sep\ \fi}
\def\printsym{\ifnum\ppsym>\z@\rsmash{^{\ppsym}}\fi}
\let\sameT=*
\def\printrcount{\count\z@=\pprcount\relax
 \ifnum\pprcount=\@ne\else
  \let\star\relax\ifnum\count\z@<\z@\count\z@-\count\z@\let\star\sameT\fi
  \rsmash{^{\the\count\z@{\star}}}%
 \fi}
\def\formsref#1{[\ref{#1}]}
\def\printformsmult{\ifnum\ppcount>\@ne\formsref{\ppforms}\let\sep\ \fi}
\newif\ifmathforms
\def\printformssingle{\ifnum\ppcount=1\ifx\ppforms\empty\else
 \ifmathforms$\m@th\ppforms$\else\ppforms\fi\let\sep\ \fi\fi}
\def\printforms{\printformsmult\printformssingle}

\def\printinfo{}
\def\finalinfo{}
\def\dobegintable{\halign\bgroup
\ \strut$\aset{##}$\hss\quad&\strut$\get@all##$&&##\cr
\noalign{\hrule\vskip2pt}
\omit\strut\quad$S$\hss&&&
 \hss$d$\hss\quad&
 \hss$n$\hss\quad&
 \ Remarks\hss\ \cr
\midstable
}
\def\doendtable{\crcr\noalign{\vskip1pt\hrule}\cr\egroup}
\def\beginstable{\hrule height\z@\hbox to\hsize\bgroup\hss\vtop\bgroup\dobegintable}
\def\endstable{\doendtable\egroup\hss\egroup}
\def\midstable{\noalign{\vskip1pt\hrule\vskip2pt}}

\def\addcomma{\ifdim\lastskip>\z@\unskip\comma\let\comma,\ \fi}
\def\see{\addcomma see\ }
\def\from#1{\addcomma from \sset{#1}}
\makeatother

\def\printinfo{&%
 \hss$\ppd$\printsym\hss\quad&
 \hss$\ppcount$\printrcount\hss\quad&%
 \printforms\printDA\sep
}
\def\finalinfo{\unskip\hss\ }

\def\step#1{\subsubsection*{Step~\rom{\iref{#1}}}}

\def\ccount{\mathrm{cc}}
\def\scount{\mathrm{sc}}
\def\ccount{c_\R}
\def\scount{s_\R}
\def\countc{_{\mathrm{c}}}
\def\counts{_{\mathrm{s}}}

\def\ts{\smash{\tilde S}}
\let\tS\ts
\def\Sym{\operatorname{Sym}}
\let\cc\Gs

\def\bh{\bold h}
\let\bh=h
\def\bts{\smash{\tilde\bS}}
\let\bts=\ts
\def\bL{\mathbf L}

\def\half{^\circ}

\def\gI{\frak I}

\makeatletter
\def\longdash{\mathrel{\dabar@\dabar@\dabar@\dabar@\dabar@\mathchar"044B}}
\makeatother

\def\barpi{\bar\pi}
\def\tp{\tilde{p}}

\let\B\barC

\let\sec\barL
\def\barU{\bar B}
\let\BU\barU

\let\BD\barD

\def\CN{B}

\def\={\mathord:\ \,}

\def\CK{\Cal K}
\def\CS{\Cal S}
\def\CM{\Cal M}
\def\ccS{\Cal S}
\def\cC{\Cal C}

\def\ttheta{\theta}

\def\tub{\operatorname{tub}}

\allowdisplaybreaks

\author{Alex Degtyarev}

\address{%
Bilkent University\\
Department of Mathematics\\
06800 Ankara, Turkey}

\thanks{%
The first author was partially supported by the
T\"{U}B\DOTaccent{I}TAK grant 
123F111.\endgraf 
The second author
was supported in part by
the ANR grants ANR-18-CE40-0009 ENUMGEOM
and ANR-22-CE40-0014 SINTROP}

\email{degt@fen.bilkent.edu.tr}

\author{Ilia Itenberg}

\address{%
Sorbonne Universit\'e and Universit\'e Paris Cit\'e, CNRS, IMJ-PRG \\
F-75005 Paris, France} 

\email{ilia.itenberg@imj-prg.fr}



\begin{document}

\begin{abstract}
We prove that the equisingular deformation type of a simple real plane sextic curve
with smooth real part is determined by its real homological type, \ie, the polarization, exceptional divisors, 
and real structure recorded in the homology of the covering $K3$-surface. As an illustration,
we obtain an equisingular deformation classification of real plane sextics with empty real part
(for completeness, we consider the few non-simple ones as well). 
\end{abstract}

\maketitle

\section{Introduction}

Recall that a {\em real algebraic variety} is a complex algebraic variety $X$
equipped with an anti-holomorphic involution $\cc\: X \to X$, referred to as the {\em real structure}
and typically omitted from the notation. The fixed point set $X_\R$ of $\cc$ is called 
the {\em real part} of $X$. If $X$ is smooth, $\dim_\C X = n$, then $X_\R$ is either empty
or a smooth manifold of real dimension $n$.
Similar terminology applies to pairs $(X, C)$, and in this paper
we deal with real (with respect to the only, up to automorphism, real structure on $\Cp2:=\Cp2_\C$,  
\viz. $[z_0 : z_1 : z_2] \mapsto [{\bar z_0} : {\bar z_1} : {\bar z_2}]$)  
algebraic curves $C \subset \Cp2$. 

From the topological point of view, the ultimate question is a deformation classification of such curves.
There is a large amount of literature on either singular {\em complex} curves or {\em smooth} real ones. 
The combination of the two, \viz. equisingular equivariant deformations of singular real curves,
is considerably less studied: worth mentioning are D. A. Gudkov \latin{et al}. \cite{Gudkov:quartics}
on irreducible quartic curves with arbitrary singularities and V. Kharlamov
\cite{Kharlamov:quintics}, where one-nodal quintics are used as a means of a deformation classification
of the smooth ones. The deformation classification of all nodal rational quintics was obtained 
by A. Jaramillo Puentes in \cite{Jaramillo.Puentes}. 

When it comes to sextics, of great help is the advanced theory of $K3$-surfaces which works
well both in singular complex and smooth real cases.
However, as was discovered in  \cite{Itenberg:LNM}, the naked notion of real homological type
(see \autoref{s.a-prohibitions} below) is no longer sufficient: in \cite{Itenberg:LNM}, 
in order to obtain the deformation classification of one-nodal real sextics,
one had to compute the fundamental polyhedra, often infinite, of certain groups generated
by reflections. This was reconfirmed by J. Josi 
\cite{Josi, Josi:arxiv}, who 
studied nodal rational real sextics; however, he observed that the problem
does not arise when none of the nodes is real. 
A generalisation of this fact to arbitrary real sextics with smooth real part
is one of the principal results of our paper, see \autoref{th.deformation-classification} below. 

From now on, to shorten the terminology, we refer to  
real curves $C \subset \Cp2$ with empty real part $C_\R = \varnothing$
as {\em empty curves}.  
The degree of such a curve (as well as that of each irreducible component thereof) 
is obviously even, and there is but one real projective equivalence class of empty conics.
The deformation class of an empty quartic is determined by its set of singularities,
which can be $\varnothing$, $2\bA_1$, $2\bA_2$, $2\bA_3$, or $4\bA_1$,
and, in the last two cases, by whether the two conic components are real or complex conjugate. 
It is also well known (by a simple convexity and codimension argument) 
that in each even degree there a unique equivariant deformation class
of smooth empty curves. 
Thus, sextics constitute the first nontrivial case.
This problem was brought to our attention by A. Libgober who was interested
in the existence of empty sextics with the set of singularities $8\bA_2$ (see \cite{Libgober:Alexander}).
We answer this question in the affirmative; moreover, we obtain a complete equisingular
equivariant deformation classification of empty sextics. 


\subsection{
Principal results}\label{s:results} 
In this paper, we mainly deal with \emph{simple} (\ie, ones with simple, \latin{aka} $\bA$--$\bD$--$\bE$, singularities)
sextics $C \subset \Cp2$ (though, see \autoref{add.J10} for the non-simple case);
it is these sextics that are closely related to $K3$-surfaces. 

One of our principal results is the equisingular equivariant deformation classification 
of empty sextics. The following theorem is proved in \autoref{proof.main}. 


\theorem\label{th.main}
There are
\roster*
\item
$169$ equisingular equivariant deformation families, contained in 
\item
$159$ real forms \rom(\latin{aka} real lattice types, see \autoref{s:RLT}\rom) of 
\item
$139$ complex lattice types \rom(see \autoref{s:complex-lattice}\rom) of 
\item
$104$ sets of singularities 
\endroster
of empty sextics\rom;
they are listed in Tables~\ref{tab.18}--\ref{tab.12}
{\rom (}see \autoref{conv.tables}{\rom )}.   
\endtheorem

\autoref{th.main} is proved lattice theoretically; 
however, in \autoref{S.trigonal} we provide an explicite geometric description of
most empty sextics with large total Milnor number. 



It is worth mentioning that, according to \cite{Aysegul:paper,degt:geography}, with the only exception 
of
the set of singularities~$2\bA_9$, $d=1$
(see \autoref{s.2A9}), each complex lattice type in 
\autoref{th.main} constitutes a single connected equisingular deformation
family. 



For completeness, in the next statement we discuss empty sextics 
with a non-simple singular point; the proof is found in \autoref{proof.J10}.

\addendum\label{add.J10}
There are two deformation families of reduced non-simple empty sextics. 
Any such sextic splits into three conics 
tangent to each other at a 
common pair
of complex conjugate points \rom(the set of singularities $2\bJ_{10}$ in the notation of~\cite{AVG1}\rom). 
One of the conics is always real, whereas the two others are either real
or complex conjugate \rom(\cf. item~\rom{\iref{2.forms}} in
\autoref{conv.tables}\rom). 
\endaddendum

\autoref{th.main} is derived from the following statement, 
proved in \autoref{s:deformation-classification},
which is of an independent 
interest. In line with the general framework of $K3$-surfaces,
we reduce the deformation classification to the purely arithmetic study
of the so-called {\em real homological types} (see \autoref{s.a-prohibitions}),  
which capture the immediate homological information about
the polarization, exceptional divisors, and real structure. 

\theorem\label{th.deformation-classification}
Two real simple sextics without real singular points are in the same
equisingular equivariant deformation class 
if and only if their real homological types
are isomorphic.
\endtheorem 


\subsection{Contents of the paper}\label{S.contents} 
In \autoref{S.prohibitions} we build the necessary algebraic/arithmetic framework
for the deformation classification of real sextics.
Upon introducing both complex and real versions of the so-called lattice and homological types,
in \autoref{s.invariants} and \autoref{s.real.invariants} we discuss their easily 
comprehensible and computable geometric invariants. 

In \autoref{S.deformations} we prove the principal results of the paper:
\autoref{th.deformation-classification} is proved first,
and \autoref{th.main} is derived therefrom by a computer-aided computation. 

In \autoref{S.trigonal}, we
describe an explicite geometric construction for the majority of special (see \autoref{ss.kernel}) 
empty sextics by means of
the double covering $p\: \Cp2 \dashrightarrow \Sigma_2$
and trigonal curves in the Hirzebruch surface~$\Sigma_2$.

In \autoref{S.explicit}, we work out a couple of examples illustrating the computation
leading to the proof of \autoref{th.main}.

\subsection{Acknowledgement}\label{S.acknowledgment}
This paper was conceived and most of its results were obtained
during our joint research stay at the {\em Max-Planck-Institut f\"ur Mathematik}, Bonn. 
We are grateful to this institution and its friendly staff
for the hospitality and excellent working conditions.

\table
\caption{The case $\mu=18$ (see \autoref{conv.tables})}\label{tab.18}
\beginstable
6A3&\relax\d{4}\conics{(3,6)}\count{1}\rcount{1}\sym{3}%
  \forms{(3,0)(0,0,3)}&\printinfo
  \see\eqref{eq.2A3+2A1.2}, \eqref{eq.6A3}
  \finalinfo\cr
2A5+4A2&\relax\d{3}\conics{(0,4)}\count{1}\rcount{1}\sym{1}%
  \forms{(1,1,1)}&\printinfo
  \see\eqref{eq.4A2}, \eqref{eq.2A5+4A2}
  \finalinfo\cr
2A5+2A3+2A1&\relax\d{2}\conics{(3,0)}\A2\count{1}\rcount{1}\sym{0}%
  \forms{(1,1)}&\printinfo
  \finalinfo\cr
2A5+2A4&\relax\d{1}\conics{}\count{1}\rcount{-2}\sym{0}%
  \forms{}&\printinfo
  \finalinfo\cr
2A6+2A2+2A1&\relax\d{1}\conics{}\count{1}\rcount{1}\sym{0}%
  \forms{}&\printinfo
  \finalinfo\cr
2A7+4A1&\relax\d{4}\conics{(3,2)}\count{1}\rcount{1}\sym{3}%
  \forms{(1,1)(1,1,0)}&\printinfo
  \see\eqref{eq.2A3+2A1.3}, \eqref{eq.2A7+4A1}
  \finalinfo\cr
2A7+2A2&\relax\d{2}\conics{(1,0)}\count{1}\rcount{1}\sym{0}%
  \forms{}&\printinfo
  \finalinfo\cr
2A8+2A1&\relax\d{1}\conics{}\count{1}\rcount{1}\sym{0}%
  \forms{}&\printinfo
  \finalinfo\cr
\same&\relax\d{3}\conics{(0,1)}\count{2}\rcount{1}\sym{1}%
  \forms{(0,1,0)-(1,0,0)}&\printinfo
  \see\eqref{eq.A8}, \eqref{eq.2A8+2A1}
  \finalinfo\cr
2A9&\relax\d{1}\conics{}\count{1}\rcount{-2}\sym{0}%
  \forms{}&\printinfo
  2 complex families
  \see\autoref{s.2A9}
  \finalinfo\cr
2D5+2A4&\relax\d{1}\conics{}\count{1}\rcount{1}\sym{0}%
  \forms{}&\printinfo
  \finalinfo\cr
2D6+2A3&\relax\d{2}\conics{(3,0)}\D2\A1\count{1}\rcount{1}\sym{1}%
  \forms{(1,1)}&\printinfo
  \see\eqref{eq.D6+2A1.2}, \eqref{eq.2D6+2A3}
  \finalinfo\cr
2D8+2A1&\relax\d{2}\conics{(3,0)}\count{1}\rcount{1}\sym{1}%
  \forms{(1,1)}&\printinfo
  \see\eqref{eq.D8.1}
  \finalinfo\cr
2D9&\relax\d{1}\conics{}\count{1}\rcount{1}\sym{0}%
  \forms{}&\printinfo
  \finalinfo\cr
2E6+2A3&\relax\d{1}\conics{}\count{1}\rcount{1}\sym{0}%
  \forms{}&\printinfo
  \finalinfo\cr
2E7+2A2&\relax\d{1}\conics{}\count{1}\rcount{1}\sym{0}%
  \forms{}&\printinfo
  \finalinfo\cr
2E8+2A1&\relax\d{1}\conics{}\count{1}\rcount{1}\sym{1}%
  \forms{}&\printinfo
  \see\autoref{s.E8}
  \finalinfo\cr
\endstable
\endtable

\table
\caption{The case $\mu=16$ (see \autoref{conv.tables})}\label{tab.16}
\beginstable
8A2&\relax\d{3}\conics{(0,4)}\count{1}\rcount{1}\sym{1}%
  \forms{(1,1,1)}&\printinfo
  \see\eqref{eq.4A2}
  \finalinfo\cr
2A3+4A2+2A1&\relax\d{1}\conics{}\count{1}\rcount{1}\sym{0}%
  \forms{}&\printinfo
  \finalinfo\cr
4A3+4A1&\relax\d{2}\conics{(3,0)}\count{1}\rcount{1}\sym{0}%
  \forms{(1,1)}&\printinfo
  \from{2D6+2A3}
  \see\autoref{ss.D6+2A1.1}
  \finalinfo\cr
\same&\relax\d{4}\conics{(3,2)}\count{4}\rcount{1}\sym{1}%
  \forms{(1,1)(1,1,0)-(3,0)(0,0,1)-(3,0)(0,2,0)-(3,0)(2,0,0)}&\printinfo
  \see\eqref{eq.2A3+2A1.1}--\eqref{eq.2A3+2A1.3}
  \finalinfo\cr
4A3+2A2&\relax\d{2}\conics{(1,0)}\count{1}\rcount{1}\sym{0}%
  \forms{}&\printinfo
  \from{6A3}
  \see\autoref{ss.2A3+2A1}
  \finalinfo\cr
2A4+2A2+4A1&\relax\d{1}\conics{}\count{1}\rcount{1}\sym{0}%
  \forms{}&\printinfo
  \finalinfo\cr
2A4+4A2&\relax\d{1}\conics{}\count{1}\rcount{1}\sym{0}%
  \forms{}&\printinfo
  \finalinfo\cr
2A4+2A3+2A1&\relax\d{1}\conics{}\count{1}\rcount{1}\sym{0}%
  \forms{}&\printinfo
  \finalinfo\cr
4A4&\relax\d{1}\conics{}\count{1}\rcount{-2}\sym{0}%
  \forms{}&\printinfo
  \finalinfo\cr
\same&\relax\d{5}\conics{(0,2)}\count{2}\rcount{1}\sym{1}%
  \forms{(0,2,0)-(2,0,0)}&\printinfo
  \see\autoref{s.2A4}
  \finalinfo\cr
2A5+6A1&\relax\d{2}\conics{(3,0)}\count{1}\rcount{1}\sym{0}%
  \forms{(1,1)}&\printinfo
  \from{2A5+2A3+2A1}
  \finalinfo\cr
2A5+2A2+2A1&\relax\d{1}\conics{}\count{1}\rcount{1}\sym{0}%
  \forms{}&\printinfo
  \finalinfo\cr
\same&\relax\d{3}\conics{(0,1)}\count{2}\rcount{1}\sym{0}%
  \forms{(0,1,0)-(1,0,0)}&\printinfo
  \from{2A5+4A2}
  \see\autoref{ss.4A2}
  \finalinfo\cr
\same&\relax\d{2}\conics{(1,0)}\count{1}\rcount{1}\sym{0}%
  \forms{}&\printinfo
  \from{2A5+2A3+2A1}
  \finalinfo\cr
\same&\relax\d{6}\conics{(1,1)}\count{2}\rcount{1}\sym{1}%
  \forms{(0,1,0)-(1,0,0)}&\printinfo
  \see\autoref{s.A5+A2+A1}
  \finalinfo\cr
2A5+2A3&\relax\d{1}\conics{}\count{1}\rcount{1}\sym{0}%
  \forms{}&\printinfo
  \finalinfo\cr
2A6+4A1&\relax\d{1}\conics{}\count{1}\rcount{1}\sym{0}%
  \forms{}&\printinfo
  \finalinfo\cr
2A6+2A2&\relax\d{1}\conics{}\count{1}\rcount{1}\sym{0}%
  \forms{}&\printinfo
  \finalinfo\cr
2A7+2A1&\relax\d{1}\conics{}\count{1}\rcount{-2}\sym{0}%
  \forms{}&\printinfo
  \finalinfo\cr
\same&\relax\d{2}\conics{(1,0)}\count{1}\rcount{1}\sym{0}%
  \forms{}&\printinfo
  \from{2D8+2A1}
  \finalinfo\cr
\same&\relax\d{4}\conics{(1,1)}\count{2}\rcount{2}\sym{1}%
  \forms{(0,1,0)-(1,0,0)}&\printinfo
  \see\autoref{s.A7+A1}
  \finalinfo\cr
2A8&\relax\d{1}\conics{}\count{1}\rcount{1}\sym{0}%
  \forms{}&\printinfo
  \finalinfo\cr
\same&\relax\d{3}\conics{(0,1)}\count{2}\rcount{1}\sym{1}%
  \forms{(0,1,0)-(1,0,0)}&\printinfo
  \see\eqref{eq.A8}
  \finalinfo\cr
2D4+2A3+2A1&\relax\d{2}\conics{(3,0)}\A1\count{1}\rcount{1}\sym{0}%
  \forms{(1,1)}&\printinfo
  \from{2D6+2A3}
  \see\autoref{ss.D6+2A1.1}
  \finalinfo\cr
2D4+2A4&\relax\d{1}\conics{}\count{1}\rcount{1}\sym{0}%
  \forms{}&\printinfo
  \finalinfo\cr
4D4&\relax\d{2}\conics{(3,0)}\count{2}\rcount{1}\sym{3}%
  \forms{(1,1)-(3,0)}&\printinfo
  \see\autoref{s.2D4}
  \finalinfo\cr
2D5+2A2+2A1&\relax\d{1}\conics{}\count{1}\rcount{1}\sym{0}%
  \forms{}&\printinfo
  \finalinfo\cr
2D5+2A3&\relax\d{1}\conics{}\count{1}\rcount{1}\sym{0}%
  \forms{}&\printinfo
  \finalinfo\cr
\same&\relax\d{2}\conics{(1,0)}\count{1}\rcount{1}\sym{0}%
  \forms{}&\printinfo
  \from{2D6+2A3}
  \finalinfo\cr
\same&\relax\d{4}\conics{(1,1)}\count{2}\rcount{1}\sym{1}%
  \forms{(0,1,0)-(1,0,0)}&\printinfo
  \see\autoref{s.D5+A3}
  \finalinfo\cr
2D6+4A1&\relax\d{2}\conics{(3,0)}\D2\count{2}\rcount{1}\sym{1}%
  \forms{(1,1)-(3,0)}&\printinfo
  \see\eqref{eq.D6+2A1.1}--\eqref{eq.D6+2A1.2}
  \finalinfo\cr
\same&\relax\d{2}\conics{(3,0)}\D1\count{1}\rcount{1}\sym{0}%
  \forms{(1,1)}&\printinfo
  \from{2D8+2A1}
  \see\autoref{ss.D8}
  \finalinfo\cr
2D6+2A2&\relax\d{1}\conics{}\count{1}\rcount{1}\sym{0}%
  \forms{}&\printinfo
  \finalinfo\cr
2D7+2A1&\relax\d{1}\conics{}\count{1}\rcount{1}\sym{0}%
  \forms{}&\printinfo
  \finalinfo\cr
2D8&\relax\d{1}\conics{}\count{1}\rcount{1}\sym{0}%
  \forms{}&\printinfo
  \finalinfo\cr
\same&\relax\d{2}\conics{(1,0)}\count{1}\rcount{2}\sym{1}%
  \forms{}&\printinfo
  \see\eqref{eq.D8.1}--\eqref{eq.D8.2}
  \finalinfo\cr
2E6+4A1&\relax\d{1}\conics{}\count{1}\rcount{1}\sym{0}%
  \forms{}&\printinfo
  \finalinfo\cr
2E6+2A2&\relax\d{1}\conics{}\count{1}\rcount{1}\sym{0}%
  \forms{}&\printinfo
  \finalinfo\cr
\same&\relax\d{3}\conics{(0,1)}\count{2}\rcount{1}\sym{1}%
  \forms{(0,1,0)-(1,0,0)}&\printinfo
  \see\autoref{s.E6+A2}
  \finalinfo\cr
2E7+2A1&\relax\d{1}\conics{}\count{1}\rcount{1}\sym{0}%
  \forms{}&\printinfo
  \finalinfo\cr
\same&\relax\d{2}\conics{(1,0)}\count{1}\rcount{1}\sym{1}%
  \forms{}&\printinfo
  \see\autoref{s.E7+A1}
  \finalinfo\cr
2E8&\relax\d{1}\conics{}\count{1}\rcount{1}\sym{1}%
  \forms{}&\printinfo
  \see\autoref{s.E8}
  \finalinfo\cr
\endstable
\endtable

\table
\caption{The case $\mu=14$ (see \autoref{conv.tables})}\label{tab.14}
\beginstable
4A2+6A1&\relax\d{1}\conics{}\count{1}\rcount{1}\sym{0}%
  \forms{}&\printinfo
  \finalinfo\cr
6A2+2A1&\relax\d{1}\conics{}\count{1}\rcount{1}\sym{0}%
  \forms{}&\printinfo
  \finalinfo\cr
\same&\relax\d{3}\conics{(0,1)}\count{2}\rcount{1}\sym{0}%
  \forms{(0,1,0)-(1,0,0)}&\printinfo
  \from{8A2}
  \see\autoref{ss.4A2}
  \finalinfo\cr
2A3+8A1&\relax\d{2}\conics{(3,0)}\A2\count{2}\rcount{1}\sym{0}%
  \forms{(1,1)-(3,0)}&\printinfo
  \from{2D6+4A1}
  \see\autoref{ss.D6+2A1.2}
  \finalinfo\cr
\same&\relax\d{2}\conics{(3,0)}\A1\count{1}\rcount{1}\sym{0}%
  \forms{(1,1)}&\printinfo
  \from{2D6+2A3}
  \see\autoref{ss.D6+2A1.1}
  \finalinfo\cr
2A3+2A2+4A1&\relax\d{1}\conics{}\count{1}\rcount{1}\sym{0}%
  \forms{}&\printinfo
  \finalinfo\cr
\same&\relax\d{2}\conics{(1,0)}\count{1}\rcount{1}\sym{0}%
  \forms{}&\printinfo
  \from{6A3}
  \see\autoref{ss.2A3+2A1}
  \finalinfo\cr
2A3+4A2&\relax\d{1}\conics{}\count{1}\rcount{1}\sym{0}%
  \forms{}&\printinfo
  \finalinfo\cr
4A3+2A1&\relax\d{1}\conics{}\count{1}\rcount{1}\sym{0}%
  \forms{}&\printinfo
  \finalinfo\cr
\same&\relax\d{2}\conics{(1,0)}\count{1}\rcount{1}\sym{0}%
  \forms{}&\printinfo
  \from{4A3+4A1}\,\eqref{eq.2A3+2A1.2}
  \see\autoref{ss.4A3+4A1}
  \finalinfo\cr
\same&\relax\d{4}\conics{(1,1)}\count{2}\rcount{1}\sym{0}%
  \forms{(0,1,0)-(1,0,0)}&\printinfo
  \from{4A3+4A1}\,\eqref{eq.2A3+2A1.1}
  \see\autoref{ss.4A3+4A1}
  \finalinfo\cr
2A4+6A1&\relax\d{1}\conics{}\count{1}\rcount{1}\sym{0}%
  \forms{}&\printinfo
  \finalinfo\cr
2A4+2A2+2A1&\relax\d{1}\conics{}\count{1}\rcount{1}\sym{0}%
  \forms{}&\printinfo
  \finalinfo\cr
2A4+2A3&\relax\d{1}\conics{}\count{1}\rcount{1}\sym{0}%
  \forms{}&\printinfo
  \finalinfo\cr
2A5+4A1&\relax\d{1}\conics{}\count{1}\rcount{1}\sym{0}%
  \forms{}&\printinfo
  \finalinfo\cr
\same&\relax\d{2}\conics{(1,0)}\count{1}\rcount{1}\sym{0}%
  \forms{}&\printinfo
  \from{2A5+2A2+2A1}
  \see\autoref{ss.A5+A2+A1}
  \finalinfo\cr
2A5+2A2&\relax\d{1}\conics{}\count{1}\rcount{1}\sym{0}%
  \forms{}&\printinfo
  \finalinfo\cr
\same&\relax\d{3}\conics{(0,1)}\count{2}\rcount{1}\sym{0}%
  \forms{(0,1,0)-(1,0,0)}&\printinfo
  \from{2A5+4A2}
  \see\autoref{ss.4A2}
  \finalinfo\cr
2A6+2A1&\relax\d{1}\conics{}\count{1}\rcount{1}\sym{0}%
  \forms{}&\printinfo
  \finalinfo\cr
2A7&\relax\d{1}\conics{}\count{1}\rcount{1}\sym{0}%
  \forms{}&\printinfo
  \finalinfo\cr
\same&\relax\d{2}\conics{(1,0)}\count{1}\rcount{2}\sym{0}%
  \forms{}&\printinfo
  \from{2A7+4A1}
  \see\autoref{ss.2A3+2A1}
  \finalinfo\cr
2D4+6A1&\relax\d{2}\conics{(3,0)}\count{2}\rcount{1}\sym{0}%
  \forms{(1,1)-(3,0)}&\printinfo
  \from{2D6+4A1}
  \see\autoref{ss.D6+2A1.2}
  \finalinfo\cr
2D4+2A2+2A1&\relax\d{1}\conics{}\count{1}\rcount{1}\sym{0}%
  \forms{}&\printinfo
  \finalinfo\cr
2D4+2A3&\relax\d{1}\conics{}\count{1}\rcount{1}\sym{0}%
  \forms{}&\printinfo
  \finalinfo\cr
\same&\relax\d{2}\conics{(1,0)}\count{1}\rcount{2}\sym{0}%
  \forms{}&\printinfo
  \from{2D6+2A3}
  \see\autoref{ss.D6+2A1.1}
  \finalinfo\cr
2D5+4A1&\relax\d{1}\conics{}\count{1}\rcount{1}\sym{0}%
  \forms{}&\printinfo
  \finalinfo\cr
\same&\relax\d{2}\conics{(1,0)}\count{1}\rcount{1}\sym{0}%
  \forms{}&\printinfo
  \from{2D6+4A1}
  \see\autoref{ss.D6+2A1.3}
  \finalinfo\cr
2D5+2A2&\relax\d{1}\conics{}\count{1}\rcount{1}\sym{0}%
  \forms{}&\printinfo
  \finalinfo\cr
2D6+2A1&\relax\d{1}\conics{}\count{1}\rcount{1}\sym{0}%
  \forms{}&\printinfo
  \finalinfo\cr
\same&\relax\d{2}\conics{(1,0)}\count{1}\rcount{1}\sym{0}%
  \forms{}&\printinfo
  \from{2D6+4A1}
  \see\autoref{ss.D6+2A1.3}
  \finalinfo\cr
2D7&\relax\d{1}\conics{}\count{1}\rcount{1}\sym{0}%
  \forms{}&\printinfo
  \finalinfo\cr
2E6+2A1&\relax\d{1}\conics{}\count{1}\rcount{1}\sym{0}%
  \forms{}&\printinfo
  \finalinfo\cr
2E7&\relax\d{1}\conics{}\count{1}\rcount{1}\sym{0}%
  \forms{}&\printinfo
  \finalinfo\cr
\endstable
\endtable

\table
\caption{The case $\mu\le12$ (see \autoref{conv.tables})}\label{tab.12}
\beginstable
12A1&\relax\d{2}\conics{(3,0)}\count{2}\rcount{1}\sym{0}%
  \forms{(1,1)-(3,0)}&\printinfo
  \see\autoref{ss.D6+2A1.2}
  \finalinfo\cr
2A2+8A1&\relax\d{1}\conics{}\count{1}\rcount{1}\sym{0}%
  \forms{}&\printinfo
  \finalinfo\cr
\same&\relax\d{2}\conics{(1,0)}\count{1}\rcount{1}\sym{0}%
  \forms{}&\printinfo
  \see\autoref{ss.D6+2A1.3}
  \finalinfo\cr
4A2+4A1&\relax\d{1}\conics{}\count{1}\rcount{1}\sym{0}%
  \forms{}&\printinfo
  \finalinfo\cr
6A2&\relax\d{1}\conics{}\count{1}\rcount{1}\sym{0}%
  \forms{}&\printinfo
  \finalinfo\cr
\same&\relax\d{3}\conics{(0,1)}\count{2}\rcount{1}\sym{0}%
  \forms{(0,1,0)-(1,0,0)}&\printinfo
  \see\autoref{ss.4A2}
  \finalinfo\cr
2A3+6A1&\relax\d{1}\conics{}\count{1}\rcount{1}\sym{0}%
  \forms{}&\printinfo
  \finalinfo\cr
\same&\relax\d{2}\conics{(1,0)}\count{1}\rcount{1}\sym{0}%
  \forms{}&\printinfo
  \see\autoref{ss.D6+2A1.1}
  \finalinfo\cr
2A3+2A2+2A1&\relax\d{1}\conics{}\count{1}\rcount{1}\sym{0}%
  \forms{}&\printinfo
  \finalinfo\cr
4A3&\relax\d{1}\conics{}\count{1}\rcount{1}\sym{0}%
  \forms{}&\printinfo
  \finalinfo\cr
\same&\relax\d{2}\conics{(1,0)}\count{1}\rcount{2}\sym{0}%
  \forms{}&\printinfo
  \see\autoref{ss.2A3+2A1}
  \finalinfo\cr
2A4+4A1&\relax\d{1}\conics{}\count{1}\rcount{1}\sym{0}%
  \forms{}&\printinfo
  \finalinfo\cr
2A4+2A2&\relax\d{1}\conics{}\count{1}\rcount{1}\sym{0}%
  \forms{}&\printinfo
  \finalinfo\cr
2A5+2A1&\relax\d{1}\conics{}\count{1}\rcount{1}\sym{0}%
  \forms{}&\printinfo
  \finalinfo\cr
\same&\relax\d{2}\conics{(1,0)}\count{1}\rcount{1}\sym{0}%
  \forms{}&\printinfo
  \see\autoref{ss.A5+A2+A1}
  \finalinfo\cr
2A6&\relax\d{1}\conics{}\count{1}\rcount{1}\sym{0}%
  \forms{}&\printinfo
  \finalinfo\cr
2D4+4A1&\relax\d{1}\conics{}\count{1}\rcount{1}\sym{0}%
  \forms{}&\printinfo
  \finalinfo\cr
\same&\relax\d{2}\conics{(1,0)}\count{1}\rcount{1}\sym{0}%
  \forms{}&\printinfo
  \see\autoref{ss.D6+2A1.3}
  \finalinfo\cr
2D4+2A2&\relax\d{1}\conics{}\count{1}\rcount{1}\sym{0}%
  \forms{}&\printinfo
  \finalinfo\cr
2D5+2A1&\relax\d{1}\conics{}\count{1}\rcount{1}\sym{0}%
  \forms{}&\printinfo
  \finalinfo\cr
2D6&\relax\d{1}\conics{}\count{1}\rcount{1}\sym{0}%
  \forms{}&\printinfo
  \finalinfo\cr
2E6&\relax\d{1}\conics{}\count{1}\rcount{1}\sym{0}%
  \forms{}&\printinfo
  \finalinfo\cr
\doendtable\egroup\quad\vtop\bgroup\dobegintable
10A1&\relax\d{1}\conics{}\count{1}\rcount{1}\sym{0}%
  \forms{}&\printinfo
  \finalinfo\cr
\same&\relax\d{2}\conics{(1,0)}\count{1}\rcount{1}\sym{0}%
  \forms{}&\printinfo
  \see\autoref{ss.D6+2A1.2}
  \finalinfo\cr
2A2+6A1&\relax\d{1}\conics{}\count{1}\rcount{1}\sym{0}%
  \forms{}&\printinfo
  \finalinfo\cr
4A2+2A1&\relax\d{1}\conics{}\count{1}\rcount{1}\sym{0}%
  \forms{}&\printinfo
  \finalinfo\cr
2A3+4A1&\relax\d{1}\conics{}\count{1}\rcount{1}\sym{0}%
  \forms{}&\printinfo
  \finalinfo\cr
\same&\relax\d{2}\conics{(1,0)}\count{1}\rcount{1}\sym{0}%
  \forms{}&\printinfo
  \see\autoref{ss.D6+2A1.1}
  \finalinfo\cr
2A3+2A2&\relax\d{1}\conics{}\count{1}\rcount{1}\sym{0}%
  \forms{}&\printinfo
  \finalinfo\cr
2A4+2A1&\relax\d{1}\conics{}\count{1}\rcount{1}\sym{0}%
  \forms{}&\printinfo
  \finalinfo\cr
2A5&\relax\d{1}\conics{}\count{1}\rcount{1}\sym{0}%
  \forms{}&\printinfo
  \finalinfo\cr
2D4+2A1&\relax\d{1}\conics{}\count{1}\rcount{1}\sym{0}%
  \forms{}&\printinfo
  \finalinfo\cr
2D5&\relax\d{1}\conics{}\count{1}\rcount{1}\sym{0}%
  \forms{}&\printinfo
  \finalinfo\cr
\midstable
8A1&\relax\d{1}\conics{}\count{1}\rcount{1}\sym{0}%
  \forms{}&\printinfo
  \finalinfo\cr
\same&\relax\d{2}\conics{(1,0)}\count{1}\rcount{1}\sym{0}%
  \forms{}&\printinfo
  \see\autoref{ss.D6+2A1.2}
  \finalinfo\cr
2A2+4A1&\relax\d{1}\conics{}\count{1}\rcount{1}\sym{0}%
  \forms{}&\printinfo
  \finalinfo\cr
4A2&\relax\d{1}\conics{}\count{1}\rcount{1}\sym{0}%
  \forms{}&\printinfo
  \finalinfo\cr
2A3+2A1&\relax\d{1}\conics{}\count{1}\rcount{1}\sym{0}%
  \forms{}&\printinfo
  \finalinfo\cr
2A4&\relax\d{1}\conics{}\count{1}\rcount{1}\sym{0}%
  \forms{}&\printinfo
  \finalinfo\cr
2D4&\relax\d{1}\conics{}\count{1}\rcount{1}\sym{0}%
  \forms{}&\printinfo
  \finalinfo\cr
\midstable
6A1&\relax\d{1}\conics{}\count{1}\rcount{1}\sym{0}%
  \forms{}&\printinfo
  \finalinfo\cr
2A2+2A1&\relax\d{1}\conics{}\count{1}\rcount{1}\sym{0}%
  \forms{}&\printinfo
  \finalinfo\cr
2A3&\relax\d{1}\conics{}\count{1}\rcount{1}\sym{0}%
  \forms{}&\printinfo
  \finalinfo\cr
\midstable
4A1&\relax\d{1}\conics{}\count{1}\rcount{1}\sym{0}%
  \forms{}&\printinfo
  \finalinfo\cr
2A2&\relax\d{1}\conics{}\count{1}\rcount{1}\sym{0}%
  \forms{}&\printinfo
  \finalinfo\cr
\midstable
2A1&\relax\d{1}\conics{}\count{1}\rcount{1}\sym{0}%
  \forms{}&\printinfo
  \finalinfo\cr
\endstable
\endtable

\section{Lattice types and homological types}\label{S.prohibitions}

In this section, we describe the complex and real versions of the so-called lattice and homological types
of a (real) simple plane sextic curve. Intuitively, at least in the complex setting, the lattice type captures
algebro-geometric properties of the sextics, whereas the homological type also takes into account
the topology of the ground field $\C$. 

In the last two subsections, we introduce the invariants of lattice/homological types
used in Tables~\ref{tab.18}--\ref{tab.12}
{\rom (}see \autoref{conv.tables}{\rom )}. 

\subsection{Complex lattice types}\label{s:complex-lattice}
Recall that, given a simple sextic $C$ in $\Cp2$,
the minimal resolution $X:=X_C$ of singularities of the double covering of $\Cp2$
ramified at $C$ is a $K3$-surface;
the covering projection is denoted by $\pi\: X \to \Cp2$.
\emph{Via} the Poincar\'{e} duality isomorphism, we always identify $H^2(X) = H_2(X)$ (unless stated otherwise,
all coefficients are in $\Z$) and regard the latter as a unimodular lattice:
$$H_2(X) \simeq \bL := 2 \bE_8 \oplus 3\bU.$$
Throughout the paper, $\bA_n$, $\bD_n$, $\bE_n$ are \emph{negative definite} root lattices
corresponding to the Dynkin diagrams of the same name, and $\bU$ is the \emph{hyperbolic plane}:
$$
\bU = \Z u_1 + \Z u_2, \qquad u_1^2 = u_2^2 = 0, \quad u_1 \cdot u_2 = 1.
$$

The classes of the exceptional divisors over the singularities of $C$ span a root lattice $S \subset H_2(X)$.
In view of the uniqueness of the decomposition of $S$ into irreducible $\bA$--$\bD$--$\bE$ components,
this lattice can be identified with the (abstract) set of singularities of $C$.
Let $h:= h_C \in H_2(X)$ be the \emph{polarization}, \ie, the class of the pull-back of a line in $\Cp2$.
One has $h^2 = 2$.
We consider the lattices
$$S_h := S \oplus \Z h,
\quad \ts_h := \text{the primitive hull of} \; S_h \; \text{in} \; H_2(X),
\quad \ts:=h^\perp\subset\ts_h.$$
The $2$-polarized hyperbolic lattice $\ts_h \ni h$ (regarded up to polarized isometry)
is called the \emph{{\rm (}complex{\rm )} lattice type} of the original sextic $C$, \cf. \cite{Shimada:splitting}.
(Here
and below,
a \emph{polarized isometry} $\Gf$
is one preserving $h$,
\ie, $\Gf(h)=h$.
If $\Gf(h)=-h$, we refer to $\Gf$ as a \emph{skew-polarized isometry}.)
If $X$ is very general in its equisingular family,
$\ts_h$ is the N\'{e}ron-Severi lattice $NS(X)$ of $X$. The Riemann-Roch theorem implies that,
in any case, $S$ is recovered
as the maximal root lattice in $h^\perp \subset NS(X)$. 

The next few statements are well known (see, \eg, \cite{degt:JAG, Shimada:splitting}) 
and follow from the global Torelli theorem \cite{Pjatecki-Shapiro.Shafarevich}
and surjectivity of the period map \cite{Kulikov:periods}
for $K3$-surfaces, as well as results of B. Saint-Donat \cite{Saint-Donat} and V. Nikulin \cite{Nikulin:Weil}; see also 
\cite{Beauville:moduli, Dolgachev:polarized}. 


\theorem\label{th.arith-conditions}
A $2$-polarized lattice $\ts_h \ni h$ represents a complex lattice type if and only if
its isomorphism class is an \emph{abstract (complex) lattice type}, \ie, 
\roster
\item $\ts_h$ is hyperbolic and admits a primitive isometry into $\bL$\rom;
\item $\ts_h$ is generated over $\Q$ by $h$ and the roots in
$h^\perp \subset \ts_h$\rom;
\item there is no element $e \in \ts_h$ such that $e^2 = 0$ and $e \cdot h = 1$.
\endroster
\endtheorem

\subsection{Real lattice types}\label{s:RLT} 
Let $C \subset \Cp2$ be a simple real sextic and $\pi\: X \to \Cp2$ the covering $K3$-surface.
The real structure $\cc$ on $\Cp2$ lifts to two commuting real structures $\cc^\pm \: X \to X$, 
so that the composition $\cc^+ \circ \cc^-$ is the deck translation $\tau$ of the covering.
The two lifts $\cc^\pm$ are distinguished by the projections $\pi(\Fix \cc^\pm) \subset \Rp2$ of their real parts;
these projections, called the \emph{halves} of $\Rp2$, are
disjoint except for the common boundary $C_\R$. 

In this paper, we consider real sextics with smooth (possibly empty) real part $C_\R$.
In this case, exactly one of the two halves of $\Rp2$ is orientable,
and the corresponding lift is denoted by $\cc^+$.
The isomorphism class of $(\ts_h \ni h, \cc_S)$, where $\cc_S \: \ts_h \to \ts_h$ is the restriction
of $\cc^+_*\: H_2(X) \to H_2(X)$,
is called the \emph{real lattice type} of
$C$.
Alternatively, the conjugacy classes, in the group of (skew-)polarized isometries of $\ts_h \ni h$,
of involutive
skew-polarized isometries $\cc_S$ as above are called the \emph{
real forms} of the complex lattice type $\ts_h \ni h$.

An obvious condition necessary for geometric realizability of a real form 
by a real sextic with smooth real part 
is that $\cc_S$
should not fix, as a set, any of the $\bA$--$\bD$--$\bE$ components of $S$.

%
%
%


\subsection{Complex homological types}\label{s:CHT}

To a simple sextic $C \subset \Cp2$ we can associate
its \emph{{\rm (}complex{\rm )} homological type}, \ie, 
the triple
$$H_2(X) \supset \ts_h \ni h$$
considered up to isometry.
An \emph{abstract {\rm (}complex{\rm )} homological type} is an isometry class of triples $\bL \supset \ts_h \ni h$, 
where $\ts_h \ni h$ is an abstract complex lattice type (see \autoref{th.arith-conditions})
and $\ts_h$
is primitive in $\bL$. 

The generic \emph{transcendental lattice} $T := \ts_h^\perp \subset \bL$ 
has two positive squares; hence, all positive definite $2$-subspaces
in $T \otimes \R$ can be oriented in a coherent way.
A choice of one of these two coherent orientations is called
an \emph{orientation} of the abstract homological type.
The homological type of a sextic has a canonical orientation, \viz.
the one given by $\R\Re\omega \oplus \R\Im\omega$, where $\omega \ne 0$
is a holomorphic $2$-form on $X$.

\theorem[see \cite{degt:JAG}]\label{th:complex-deformation}
The equisingular deformation families of simple sextics
are in a canonical bijection with
the oriented abstract homological types.
\endtheorem


\subsection{Real homological types}\label{s.a-prohibitions}
The \emph{real homological type} of a real sextic $C \subset \Cp2$
with smooth real part
is the quadruple
\[
\bigl(H_2(X) \supset \ts_h \ni h, \cc^+_*\bigr)\label{real-hom-type}
\]
considered up to isometry commuting with $\cc^+_*$.
Here, the orientation of the complex homological type is redundant,
as $-\cc^+_*$ is an orientation reversing isometry.

In this paper, we are particularly interested in empty sextics.

\theorem\label{th.real-statement}
A quadruple $(\bL \supset \ts_h \ni h,
\ttheta)$ represents the real homological type
of an empty sextic if and only if
\roster
\item
$\bL \supset \ts_h \ni h$ is an abstract complex homological type,
\item
$\ttheta$ is an involutive skew-polarized isometry of $\bL \supset \ts_h \ni h$, and
\item\label{i.Nikulin} $\Ker(1 - \ttheta) \simeq \bE_8(2) \oplus \bU(2)$.
\endroster
\endtheorem

\proof
The first two requirements are obvious, whereas the necessity of item \iref{i.Nikulin}
follows from the deformation classification of smooth real sextics found in \cite{Nikulin:forms}.
Conversely, in view of \autoref{th:complex-deformation}, the sufficiency is given by the construction
of anti-holomorphic maps (see, \eg, \cite{Barth.Peters.VanDeVen} or \cite[\S 13.4.3]{DIK}
for an explanation in the anti-holomorphic setup);
since the real structure obtained preserves the polarisation,
it automatically commutes with the deck translation.
\endproof

\definition\label{def.empty.sphere}
If the real sextic $C \subset \Cp2$ is empty,
there is a distinguished sphere $\Fix\cc^- \subset X$ disjoint from all exceptional divisors. 
Choosing an orientation, we obtain  
a class $s:=[\Fix\cc^-] \in H_2(X)$, called the \emph{empty sphere}.
It is immediate that $s$ is $\cc^+_*$-skew-invariant,
and
\[*
s^2=-2,\quad s\cdot \ts_h = 0,\quad s=h\bmod 2 H_2(X);
\]
since $h^2=2$, the primitive hull~$U_{h, s}$ of $\Z h\oplus\Z s$ in~$H_2(X)$ is isomorphic
to~$\bU$ and, thus, splits as an orthogonal direct summand of $H_2(X)$,
orthogonal to $\ts = h^\perp\subset\ts_h$.
\enddefinition

\corollary\label{cor.reduction}
Let $\ts_h \ni h$ be a complex lattice type and $\ttheta_S\: \ts_h \to \ts_h$
an involutive 
skew-polarized isometry, and let $\ttheta_h$ be the restriction of~$\ttheta_S$ to
$\ts$.
If
$(\ts_h \ni h, \ttheta_S)$ represents the real lattice type of an empty
sextic,
then
\begin{align}\label{i.embeddable}
& \text{both
$\ts_\pm:=\Ker(1\mp\ttheta_h)$
admit primitive isometries into $\bE_8(2)\oplus\bU(2)$,} \\
\label{i.orthogonal}
& \text{$\Z h$ is an orthogonal direct summand in $\ts_h$;}
\end{align}
here and below, given a lattice $M:=(M, b)$, where $b$ is the bilinear form on $M$, 
the notation $M(2)$ stands for the lattice $(M, 2b)$. 
\endcorollary

\proof
If
$\ttheta_S$ extends to an empty real structure $\ttheta = \cc^+_*$,
by \autoref{th.real-statement} we have
\[\label{L+L-}
\aligned
L_+&:=\Ker(1 - \cc^+_*) \simeq \bE_8(2) \oplus \bU(2),\\
L_-&:=\Ker(1 + \cc^+_*) \simeq \bE_8(2) \oplus \bU(2) \oplus\bU,
\endaligned
\]
and statement~\eqref{i.embeddable} is immediate for
the primitive sublattice $\ts_+\subset L_+$.
For $\ts_-$, we observe in addition that it is orthogonal to
$U_{h, s} \subset L_-$,
and
$U_{h, s}^\perp \subset L_-$ equals $\bE_8(2) \oplus \bU(2)$ 
in view of the uniqueness of the latter lattice in its genus.
\endproof

\remark\label{rem.double}
As an immediate consequence of \autoref{cor.reduction}, we have
$\ts_\pm=\ts\half_\pm(2)$ for some even lattices
$\ts\half_\pm$, and~\eqref{i.embeddable}
reduces to the existence of primitive isometries
$\ts\half_\pm \into \bE_8\oplus\bU$, which is
easily checked using Nikulin's theory~\cite{Nikulin:forms}.
\latin{A posteriori},
our \autoref{th.main} implies
that the necessary condition \eqref{i.embeddable}
is also sufficient,
provided that the existence of the complex family is known. 
In the case of non-special sextics (see \autoref{ss.kernel}),
this fact has a simple direct proof, see \autoref{th.non-special} below.
\endremark 

\subsection{Invariants of complex lattice types}\label{s.invariants}
Fix a complex lattice type $\ts_h\ni h$ and pick a representative sextic
$C\subset\Cp2$. The following
objects are invariants of $\ts_h\ni h$. 

\subsubsection{\rm The \emph{kernel} of the extension,}\label{ss.kernel}
\ie, the (finite) isotropic subgroup
\[*
\ts_h/S_h\subset\discr S_h:=S_h^\vee\!/S_h.
\]
In a sense, given~$S$, this kernel \emph{is} the lattice type; however, we
confine ourselves to the exponent
(maximal order of elements)
\[*
d:=\exp(\ts_h/S_h)
\]
and use a number of
more geometric/numeric derivatives described below.
The sextic is reducible if and only if $2\mathrel|d$, see \cite{degt:JAG}. 
Sextics with $3\mathrel|d$ are often said to be of \emph{torus type}.
In general, a sextic is called \emph{special} if $d>1$.

\subsubsection{\rm The \emph{combinatorial type},}\label{ss.combinatorial}
\ie, the homeomorphism class of the pair $(\tub_C,C)$, where $\tub_C$ is a regular
tubular neighbourhood of~$C$. Roughly, 
this consists of (the degrees of) the
irreducible components of~$C$, its set of singularities~$S$, and the
position of (the branches of) the singular points of~$C$ on its components.

If $C_\R$ is smooth, all components of~$C$ must be of even degree; hence, these
degrees $(6)$, $(4,2)$, or $(2,2,2)$ are determined by the number $c_\C$ of
conic components.

\subsubsection{\rm The \emph{alignment},}\label{ss.alignment}
which is a partial description of
the ``position of the singular points on the components''
mentioned in the previous section.
Assume
that $c_\C=3$ and consider a singular point $P$ of type
$\bA_{2p+1}$ or $\bD_{2p+4}$, $p\ge1$. Then, two of the three conic
components of~$C$ are tangent to each other at~$P$; they are called
\emph{tight} at~$P$, whereas the third conic is called \emph{loose}.
Now, assume that $C$ has exactly two points $P_1,P_2$ of
the same type $\bA_{2p+1}$ or $\bD_{2p+4}$, $p\ge1$. 
If $p=1$, there are two possibilities,
\viz.
\roster
\def\theenumi{\roman{enumi}}
\item\label{DA-1}
the points are \emph{misaligned}, \ie,
the loose conics at $P_1,P_2$ are distinct, or
\item\label{DA-2}
the points are \emph{aligned}, \ie,
$P_1,P_2$ have a common loose conic,
\endroster
which are usually indicated \via, \eg, \str{\bA_3} or \sstr{\bD_6}.
If $p>1$, the pair is automatically misaligned, as in~(\ref{DA-1}),
and this property
is obviously preserved under perturbations (provided that the new curve still splits into
three conics). 

If $c_\C=1$, one can also consider the set of singularities of the quartic
component, but we use this extra piece of data only once,
in \autoref{rem.aposteriori} below. 

\subsubsection{\rm Splitting conics and lines}\label{ss.splitting}
A conic or line $\CN\subset\Cp2$ is \emph{splitting}
(\emph{$Z$-splitting} in~\cite{Shimada:splitting}) for~$C$ if
\roster*
\item
the pull-back $\pi\1(\CN)$ splits into two smooth rational curves
$\CN'$, $\CN''$, and
\item
the classes $[\CN']$, $[\CN'']$ are distinct and lie in $\ts_h\subset H_2(X)$.
\endroster
The latter condition ensures that the splitting conics and lines are stable
in the sense that they follow all equisingular deformations of~$C$.

\theorem[Shimada~\cite{Shimada:splitting}]
A complex lattice type is determined by its combinatorial type and the numbers of
splitting lines and conics. 
\endtheorem

In view of \eqref{i.orthogonal}, an empty sextic cannot have splitting lines;
therefore,
it suffices to consider the number $s_\C$ of its splitting conics.

\remark\label{rem.conics}
Whenever \eqref{i.orthogonal} holds,
both conic components and splitting conics are found using vectors
$v\in\tS/S$ whose shortest representative in~$\tS$ has square~$-4$.
The conic components are in a bijection with such vectors of order~$2$
(equivalently, invariant under the deck translation), whereas splitting
conics correspond to pairs of opposite vectors of any larger order
(equivalently, $2$-element
orbits of the deck translation); in the latter case, the individual vectors 
are in a bijection with the pull-backs of the splitting conics in~$X$.
Splitting conics of order~$3$ are often referred to as
\emph{torus structures};
if $\{f=0\}$ is such a conic, the equation of~$C$
is
$f^3+g^2=0$
for some cubic polynomial $g$. 

As a consequence, the presence of a splitting conic implies that $d>2$.
\endremark


\subsubsection{\rm The group $\Sym(\ts_h\ni h)$ of \emph{stable symmetries},}\label{ss.sym}
\ie, the subgroup of $\OG_h(\ts_h)$ preserving the exceptional divisors (as a
set) and acting identically on $\discr\ts_h$. If $C$ is generic,
so that $\ts_h$ is the N\'{e}ron-Severi lattice of $X$,
this is indeed the group of symplectic automorphisms of~$X$
commuting with the deck translation. These automorphisms descend to~$\Cp2$
and are stable in the sense that they follow all equisingular
deformations of~$C$; we refer to~\cite{degt:symmetric} for further details.

We are mainly interested in nontrivial \emph{stable involutions};
the number of such involutions is
denoted by~$m_\C$.

\subsection{Invariants of real forms}\label{s.real.invariants}
Given a real form~$\cc_S$ 
of a complex lattice type $\ts_h \ni h$, the invariants
introduced in \autoref{s.invariants} have the following real refinements. 

The conic count~$c_\C$ splits into $\ccount:=(r\countc,c\countc)$, where
\roster*
\item
$r\countc$ is the number of real conic components and
\item
$c\countc$ is the number of pairs of complex conjugate ones,
\endroster
so that $c_\C=r\countc+2c\countc$.
Similarly, we let $\scount:=(r\counts,c\counts,q\counts)$, where
\roster*
\item
$r\counts$ is the number of real splitting conics whose pull-backs in~$X$ are
also real,
\item
$c\counts$ is the number of real splitting conics whose pull-backs in~$X$ are
complex conjugate to each other, and
\item
$q\counts$ is the number of pairs of complex conjugate splitting conics,
\endroster
so that $s_\C=r\counts+c\counts+2q\counts$.

Besides, we have the group $\Sym_\R(\ts\ni h,\cc_S)$ of equivariant stable
symmetries and the number $m_\R$ of equivariant nontrivial stable
involutions, \cf. \autoref{ss.sym}.

Now, we are ready to describe the data presented in Tables~\ref{tab.18}--\ref{tab.12}.

\convention\label{conv.tables}
The rows of Tables~\ref{tab.18}--\ref{tab.12} list
all complex
lattice types $\ts_h\ni h$, $S\ne0$, admitting at least one real form represented 
by an empty sextic. 

The $S$-column refers to the set of singularities~$S$. 

The $d$-column shows the exponent $d=\exp(\ts_h/S_h)$ (see
\autoref{ss.kernel}) and the number $m_\C=m_\R$ (if greater than~$0$) of
stable involutions (see \autoref{ss.sym}) as a superscript.

The $n$-column shows the number~$n$ of real forms of
$\ts_h\ni h$
and the
number of deformation families (if greater than~$1$) \emph{per real form}
as a superscript; the latter is followed by a $\sameT$ if all 
real homological types
(within one real form)
share the same
equivariant transcendental lattice
$(T, \cc_T)$.

Finally, in the remark column, we describe the invariants of the real forms:
\roster*
\item
$\ccount$ (if $n=1$ and $c_\C=3$), as a two element list $(*,*)$,
\item
$\scount$ (if $n=1$ and $s_\C>0$), as a three element list $(*,*,*)$, and
\item
the alignment (see \autoref{ss.alignment}), whenever applicable.
\endroster
(To save space, we omit the labels $\ccount$, $\scount$.)
If $n>1$, we encounter but
the four cases below,
and the multiple values of~$\ccount$ and~$\scount$
are replaced with a reference
\formsref{2.forms}--\formsref{4.forms}:
\roster
\item\label{(1,1)-(3,0)}\label{2.forms}
one has $s_\C=0$ and
the two real forms differ by $\ccount=(1,1)$ or $(3,0)$;
\item\label{(0,1,0)-(1,0,0)}
one has $c_\C\le1$ and
the two real forms differ by $\scount=(0,1,0)$ or $(1,0,0)$;
\item\label{(0,2,0)-(2,0,0)}
one has $c_\C=0$ and
the two real forms differ by $\scount=(0,2,0)$ or $(2,0,0)$;
\item\label{(1,1)(1,1,0)-(3,0)(0,0,1)-(3,0)(0,2,0)-(3,0)(2,0,0)}\label{4.forms}
there are four real forms, one with $\ccount=(1,1)$ and $\scount=(1,1,0)$
and three with
$\ccount=(3,0)$ and $\scount=(0,0,1)$, $(0,2,0)$, or $(2,0,0)$.
\endroster

In the same column, we explain the construction of real curves, by a
perturbation (see
\autoref{rem.perturbations} below) from a larger set of
singularities $S'$ (\via\ ``from $S'$'') and/or by a reference to \autoref{S.trigonal} or
\autoref{S.explicit}. The common reference for most non-special ($d=1$)
curves is \autoref{s.d=1}.
\endconvention 

\remark\label{rem.perturbations}
Another well-known consequence of the standard theory of $K3$-surfaces 
 \cite{Pjatecki-Shapiro.Shafarevich, Kulikov:periods, Saint-Donat} 
is the fact that perturbations of simple sextics are unobstructed. 
Literally the same argument shows that the statement holds
in the equivariant setting as well. More precisely, let \eqref{real-hom-type}
be the real homological type of a real sextic with smooth real part.
(The last condition is not essential, but then one would have to make a choice between $\sigma^\pm$.) 
Pick a primitive $\sigma^+_*$-invariant root sublattice $S' \subset S$, 
and denote by $\ts'_h$ the primitive hull of $S' \oplus \Z h \subset \ts_h$.
Then, there is a family of real sextics $C_t$, $t \in [0, 1]$, 
such that $C_0 = C$ and the real homological type of each $C_t$, $t \in (0, 1]$,  
is $\bigl(H_2(X) \supset \ts'_h \ni h, \cc^+_*\bigr)$. 
\endremark

\remark\label{rem.aposteriori}
\latin{A posteriori}, we conclude that the chosen invariants do completely
describe our results.

First, whenever a complex lattice type has several real forms, these forms
are distinguished by the pairs $(\ccount,\scount)$, see
\iref{2.forms}--\iref{4.forms} above.

Second, for the sets of singularities \emph{listed in the tables} and complex lattice 
types without components of odd degree, the invariants
$(d,c_\C,s_\C)$ and the alignment~\iref{DA-1}
\vs. \iref{DA-2} in \autoref{ss.alignment}
(whenever applicable)
almost single out those admitting at least one empty real
structure. The two exceptions are
the sets of singularities
\sset{4A3+2A1}
and \sset{2A3+6A1}, both with $c_\C=1$: in addition, we obviously need to
exclude the complex lattice types where the quartic component has a single type 
$\bA_3$ point.

Finally, we observe that, for all empty sextics, one has $\Sym_\R=\Sym$.
\endremark

\section{Deformation classification}\label{S.deformations}

In general, it is not true (\cf.~\cite{Itenberg:LNM})
that the equisingular equivariant deformation class 
of a real sextic $C \subset \Cp2$ is determined by its real homological type.
It is, however, true in the special case where $C$ has no real singular points
(\cf. \cite{Josi:arxiv} for nodal sextics); this fact and its implications for
empty real sextics are the principal results of this section.

\subsection{Proof of Theorem \ref{th.deformation-classification}}\label{s:deformation-classification}



Fix a quadruple
\[
(\bL \supset \bts_\bh \ni \bh, \ttheta)
\label{eq.fixed}
\]
and assume that
it represents the real homological type of a real sextic without real singular points.
Fix
a $(-\ttheta)$-invariant
Weyl chamber $\Delta$ of $\bh^\perp \subset \bts_\bh$,
which we regard as a distinguished set of vectors in $\bh^\perp$.
Denote by $S:=\Z\Delta$ the sublattice generated by the roots in $\bh^\perp$.

Put $L_\pm = \Ker(1 \mp \ttheta)$ and consider the
\roster*
\item hyperbolic lattices $P_\pm := (\bts_\bh \cap L_\pm)^\perp \subset L_\pm$,
\item projectivized positive cones
$\cC_\pm := \bigl\{x \in P_\pm \otimes \R \bigm | x^2 = 1\bigr\}$, and
\item walls
$H_v := \bigr\{(x_+,x_-)\in\cC_+\times\cC_-\bigm|x_+\cdot v =x_-\cdot v=0\bigl\}$,
where $v \in \bL$ is a vector such that $v^2 = -2$ and $v \cdot \bh = 0$. 
\endroster

A \emph{marking} of a real sextic curve $C \subset \Cp2$ as above
is an isometry $\psi\: H_2(X_C) \to \bL$ such that
\roster*
\item $\psi$ is real, \ie, $\psi \circ \cc^+_* = \ttheta \circ \psi$,
\item $\psi(h_C) = \bh$, and
\item$\psi$ establishes a bijection between the exceptional divisors of $X_C$
and $\Delta$.
\endroster
The \emph{period space} of marked real sextics is the space
$$
\Omega := \bigl(\cC_+ \times \cC_- \sminus \textstyle{\bigcup H_v}\bigr)/ {\pm1},
$$
and the \emph{period map} sends a sextic $C$ to the class of $\psi(\omega_+, \omega_-)$,
where $\omega_+ = \Re\omega$, $\omega_-= \Im\omega$,
and $\omega$ is a non-zero real (in the sense that $\cc^+_*(\omega) = \bar\omega$) holomorphic $2$-form on $X_C$.
The period map makes $\Omega$ a fine moduli space of marked real simple sextics
with the given real homological type~\eqref{eq.fixed}:
in view of \cite{Beauville:moduli} and \cite{Saint-Donat}, this fact follows from the construction
of anti-holomorphic maps, as explained in the proof of \autoref{th.real-statement}.
Certainly, this space is disconnected,
and the rest of the proof deals with showing that the connected components of $\Omega$
constitute a single orbit of the automorphism group of~\eqref{eq.fixed};
in other words, any component can be sent
to any other by a change of the marking.

First, observe that each of $\cC_\pm$ has two connected components. Hence, so
does $(\cC_+\times\cC_-)/{\pm1}$, and these two components are interchanged
by~$-\ttheta$.


Now, consider a wall $H_v$ for some $v\in\bL\sminus S$, $v^2=-2$,
$v\cdot\bh=0$. It is obvious that $\codim H_v\ge2$ unless $v \in L_\pm$, and
$H_v=\varnothing$ unless
\[
\text{the lattice $S+\Z v$ is negative definite}.
\label{eq.neg.def}
\]
Therefore, from now on we assume that $v \in L_\pm$ and \eqref{eq.neg.def}
holds. In particular, the latter implies that $\ls|v\cdot e|\le1$ for each
$e\in\Delta$. To complete the proof, we need to find,
for any such vector~$v$, an automorphism of~\eqref{eq.fixed}
that would interchange
pairs of
components adjacent to~$H_v$ (more precisely, any
pair of
components whose
closures share a common
very general
point of~$H_v$).



If $v \cdot e = 0$ for each vector $e \in \Delta$,
the reflection
\[
r_v\:x \mapsto x + (x \cdot v)v
\label{eq.reflection}
\]
is the desired automorphism.

In general, we depict the union $\Delta\cup v$ by an analog~$D_v$ of Dynkin
diagram in which dotted edges $[e_1,e_2]$ (corresponding to the value
$e_1\cdot e_2=-1$) are allowed. To minimize the number of such edges, we
occasionally change the signs of the vectors. In particular, assuming that
$v \in L_\pm$, we let $\bar{v}:=v$ and $\bar{e}:=\mp\ttheta(e)$ for $e \in S$,
so that
$[\bar{e}_1,\bar{e}_2]$ is always an edge of the same type
(solid, empty, or dotted) as $[e_1,e_2]$.
Due to the assumption that $C$ has no real singular points, the Dynkin diagram of
$\Delta$ itself is a disjoint union of trees
(of types $\bA$--$\bD$--$\bE$) split into pairs exchanged by
the complex conjugation; thus, changing, if necessary (\ie, if $v \in L_+$),
the signs in half of the components, we can assume that these pairs are of the
form $\Sigma,\bar{\Sigma}$.

\figure
\begin{picture}(51,54)(-10,-14)
\def\c{\circle4}
\def\b{\circle*4}
\def\h{\line(1,0){10}}
\def\v{\line(0,1){10}}
\let\e=e
\def\be{\bar e}
\def\clap#1{\hbox to0pt{\hss$#1$\hss}}
\put(0,0)\c
\put(0,-12){\clap{\be_{n}}}
\put(2,0)\h
\put(14,0)\c
\put(18,-1){\hbox to18pt{\hss\dots\hss}}
\put(38,0)\c
\put(38,-12){\clap{\be_1}}
\put(0,2)\v
\put(0,14)\b
\put(0,16)\v
\put(0,28)\c
\put(0,35){\clap{\e_{n}}}
\put(2,28)\h
\put(14,28)\c
\put(18,27){\hbox to18pt{\hss\dots\hss}}
\put(38,28)\c
\put(38,35){\clap{\e_1}}
\end{picture}\qquad
\begin{picture}(40,54)(-23,-14)
\def\clap#1{\hbox to0pt{\hss$#1$\hss}}
\put(0,24){\clap{\Sigma\simeq\bA_{n}}}
\put(0,-2){\clap{\bar\Sigma=\mp\ttheta(\Sigma)}}
\end{picture}
\qquad
\begin{picture}(75,54)(-10,-14)
\def\c{\circle4}
\def\b{\circle*4}
\def\h{\line(1,0){10}}
\def\v{\line(0,1){10}}
\let\e=e
\def\be{\bar e}
\def\clap#1{\hbox to0pt{\hss$#1$\hss}}
\def\l{\clap{\vtop to10pt{\parindent0pt\leaders
 \vbox{\vskip.5pt\vrule height1.5pt \vskip.5pt}\vfil}}}
\put(0,0)\c
\put(0,35){\clap{\e_{n}}}
\put(4,-1){\hbox to18pt{\hss\dots\hss}}
\put(24,0)\c
\put(24,35){\clap{\e_k}}
\put(26,0)\h
\put(38,0)\c
\put(42,-1){\hbox to18pt{\hss\dots\hss}}
\put(62,0)\c
\put(62,35){\clap{\e_1}}
\put(24,12)\l
\put(25,15){\line(1,1){11.5}}
\put(24,14)\b
\put(25,13){\line(1,-1){11.5}}
\put(24,26)\l
\put(0,28)\c
\put(0,-12){\clap{\be_{n}}}
\put(4,27){\hbox to18pt{\hss\dots\hss}}
\put(24,28)\c
\put(24,-12){\clap{\be_k}}
\put(26,28)\h
\put(38,28)\c
\put(42,27){\hbox to18pt{\hss\dots\hss}}
\put(62,28)\c
\put(62,-12){\clap{\be_1}}
\end{picture}
\caption{Extended graphs $D_v$}\label{fig.graphs}
\endfigure

Observe that, in any induced \emph{subtree} $T\subset D_v$, the signs
can be changed so that $T$ has no dotted edges; then, by~\eqref{eq.neg.def},
$T$ must be an ordinary simply laced Dynkin diagram.
We conclude that $v$ is adjacent to at most one
pair of components $\Sigma,\bar{\Sigma}$, as otherwise $D_v$ would contain $\tD_4$.
Likewise, if $v$ is adjacent to a pair of vertices $e_1,e_2\in\Sigma$
(hence, also to $\bar{e}_1,\bar{e}_2\in\bar\Sigma$), then
(\cf. \autoref{fig.graphs}, right)
\roster*
\item
$e_1$ and $e_2$ are adjacent in $\Sigma$, as otherwise $D_v\supset\tD_4$,
\item
$[v,e_1]$ and $[v,e_2]$ are edges of the opposite types, as otherwise
$D_v\supset\tA_2$,
\item
as a consequence, $v$ is not adjacent to any other vertex $e\in\Sigma$.
\endroster
Finally, we have $\Sigma\simeq\bA_n$, $n\ge1$, as otherwise $D_v$ would
contain a subgraph $\tD_m$. Summarizing, we arrive at the two configurations
shown in \autoref{fig.graphs}. In the left figure, the union
$\Sigma\cup v\cup\bar\Sigma$ is $\bA_{2n+1}$, with the standard basis
\[*
e_1,\ldots,e_n,e_{n+1}:=v,e_{n+1}:=\bar{e}_n,\ldots,e_{2n+1}:=\bar{e}_1.
\]
In the right figure, the \emph{lattice} $\Z\Sigma+\Z v+\Z\bar\Sigma$ is
still~$\bA_{2n+1}$: we merely change the definition of
\[*
e_{n+1}:=v-(e_k+\ldots+e_n+\bar{e}_n+\ldots+\bar{e}_k).
\]
In both cases, we need an element~$r$
of the Weyl group of $\bA_{2n+1}$ commuting with~$\ttheta$,
preserving (as a set) the subset
$\{e_1,\ldots,e_n,\mp e_{n+1},\ldots,\mp e_{2n+1}\}$, and acting \via\ $-1$ on its
orthogonal complement. If $v\in\bL_+$ (the ``$-$'' sign above), then
$-r=\ttheta$
is the automorphism induced by the only nontrivial symmetry of the Dynkin
graph~$\bA_{2n+1}$.
If $v\in\bL_-$ (the ``$+$'' sign), then $r$ is
\[*
e_i\leftrightarrow e_{i+n+1},\ i=1,\ldots,n,\qquad
e_{n+1}\mapsto-(e_1+\ldots+e_{2n+1}).
\]
Instead of decomposing~$r$ into a product of reflections, we
observe
that
it preserves $\sum ie_i\bmod(2n+2)\bA_{2n+1}$. Thus, $r=\id$
on $\discr\bA_{2n+1}$ and, therefore, $r$ extends
to~$\bL$ identically on the
orthogonal complement $\bA_{2n+1}^\perp$.
\qed

\subsection{Proof of \autoref{th.main}}\label{proof.main}
\autoref{th.deformation-classification} reduces the proof to the enumeration of
appropriate real homological types, which is done in several steps:
\roster
\item\label{step.complex}
we list all
triples $(\ts_h \ni h, \ttheta_S)$
satisfying
\iref{i.embeddable} and \iref{i.orthogonal};
\item\label{item.T}
for each triple, we list all
eigenlattices $T_\pm$ of the prospective equivariant transcendental lattices $(T, \ttheta_T)$;
\item\label{item.cc.T}
for each pair $(T_+, T_-)$,
the prospective transcendental lattices $(T, \ttheta_T)$ are obtained as appropriate
finite index extensions of $T_+ \oplus T_-$;
\item\label{item.cosets}
we list the isomorphism classes of unimodular equivariant finite index extensions $(\bL, \ttheta)$
of $(\ts_h \oplus T, \ttheta_S \oplus \ttheta_T)$;
\item\label{item.sanity.check}
from the real homological types $(\bL \supset \ts_h \ni h, \ttheta)$
thus obtained
we select those representing empty sextics, see \autoref{th.real-statement}\iref{i.Nikulin}.
\endroster

\step{step.complex}
We start with the known classification of complex lattice types of simple sextics \cite{Yang}
and select those satisfying obvious
restriction, \viz. the fact that
\[
\text{the number of singular points of each type is even}.
\label{eq.S.even}
\]

For each complex lattice type $\ts_h$ thus obtained,
we list the conjugacy classes of skew-polarized involutions
$\ttheta_S\: \ts_h \to \ts_h$
such that
\[
\text{$\ttheta_S$ does not fix as a set any irreducible summand of $S$}.
\label{eq.no.fixed}
\]
Then, we select those pairs $(\ts_h, \ttheta_S)$ that satisfy \eqref{i.embeddable},
see \autoref{rem.double}.

\remark
\latin{A posteriori}, we confirm 
that any complex lattice type $\ts_h\ni h$
that admits an involution as in~\eqref{eq.no.fixed} has $\Z h$ as an
orthogonal direct summand, \cf.~\eqref{i.orthogonal};
hence, there are neither linear/cubic 
components nor splitting lines. 

Furthermore,
within each set of singularities,
any two complex lattice types that admit an involution as in \eqref{eq.no.fixed} are 
distinguished by the invariants $(d,c_\C,s_\C)$ and alignment
(see \autoref{s.invariants}), whereas within
each complex lattice type, the involutions satisfying~\eqref{eq.no.fixed}
are distinguished by the invariants $(c_\R,s_\R)$
introduced in \autoref{s.real.invariants}. In this sense, our tables are
complete. 

Moreover, with two exceptions, the complex lattice types admitting an
involution as in~\eqref{eq.no.fixed} are distinguished by
$(d,c_\C,s_\C)$ and alignment from those not admitting one. The exceptions
are the sets of singularities \singset{4A3+2A1} and \singset{2A3+6A1} with a
single conic component passing through an \emph{odd} number of $\bA_3$-type
points (or, in other words, a quartic component with a single $\bA_3$-type
singularity).
\endremark

\remark\label{rem.always.double}
Another experimental observation is the fact that, for any involution
satisfying~\eqref{eq.no.fixed}, the eigenlattices $\ts_\pm$ are
of the form $\ts\half_\pm(2)$ (\cf. \autoref{rem.double}), where
$\ts\half_\pm$ are \emph{even root lattices}.
\endremark

\step{item.T}
Let $\ts_\pm = \ts\half_\pm(2)$, see \autoref{rem.double} or \autoref{rem.always.double},
pick a primitive isometry $\ts\half_\pm \into \bE_8 \oplus \bU$, and denote by $T\half_\pm$
the orthogonal complement.
Then, we have $T_+ = T\half_+(2)$ and $T_- = T\half_-(2) \oplus \Z s$, $s^2= -2$,
see \autoref{def.empty.sphere}.
This construction determines the genera of $T_\pm$ (see, \eg, \cite{Nikulin:forms}),
and we check, on a case by case basis, that each lattice obtained is unique in its genus.

Since it is easier to work with lattices rather than their discriminant forms,
we construct an isometry $\ts\half_\pm \into \bE_8 \oplus \bU$ by
representing the Dynkin diagram of the root system of $\ts\half_\pm$,
see \autoref{rem.always.double}, as an induced subgraph
of the graph
\[
\vcenter{\hbox{\begin{picture}(116,22)(-2,-4)
\def\c{\circle4}
\def\b{\circle*4}
\def\h{\line(1,0){10}}
\def\v{\line(0,1){10}}
\def\clap#1{\hbox to0pt{\hss$#1$\hss}}
\put(0,14)\b
\put(2,14)\h
\put(14,14)\b
\put(16,14)\h
\put(28,14)\b
\put(30,14)\h
\put(28,0)\b
\put(28,2)\v
\put(42,14)\b
\put(44,14)\h
\put(56,14)\b
\put(58,14)\h
\put(70,14)\b
\put(72,14)\h
\put(84,14)\b
\put(86,14)\h
\put(98,14)\b
\put(100,14)\h
\put(112,14)\b
\end{picture}}}
\label{eq.E10}
\]
which is the Coxeter scheme of a fundamental polyhedron
of the group generated by reflections of $\bE_8 \oplus \bU$,
see, \eg, \cite{Vinberg:polyhedron}.

For the uniqueness in the genus, we encounter three cases:
\roster*
\item definite lattices of rank $1$: the uniqueness is obvious;
\item hyperbolic lattices of rank $2$: we are not aware
of any general theory, but our needs are completely covered
by the reduced form found, \eg, in~\cite{Dickson};
\item hyperbolic lattices of higher rank:
we use Miranda--Morrison's theory \cite{Miranda.Morrison:book}.
\endroster

\step{item.cc.T} Here and at step \iref{item.cosets}, we need to solve the following problem:
given two non-degenerate even lattices $M_1$ and $M_2$, find all (up to $O(M_1) \times O(M_2)$ or a certain
prescribed subgroup thereof) classes of even finite index extensions 
\[\label{eq.primitive-extension} 
\text{${\tilde M} \supset M_1 \oplus M_2$ such that both $M_1$ and $M_2$
are primitive in ${\tilde M}$.}   
\]
To this end, recall that a non-degenerate lattice $M$ is naturally a subgroup of 
its dual lattice $M^\vee \subset M \otimes \Q$, and the \emph{discriminant form} of a non-degenerate even lattice $M$
is the finite abelian group $\discr M := M^\vee \! / M$ equipped with the induced $\Q/2\Z$-valued
quadratic form
(see \cite{Nikulin:forms}, where the notation is $q_M$).
We denote by $\gI_*(M)$ the image of the natural homomorphism
$\OG_*(M) \to \Aut_*(\discr M)$, where $*$ is a placeholder for a number of extra symbols
used below to restrict the groups, \eg, preserving $h$, commuting with an involution, \etc.

According to \cite[Proposition 1.4.1]{Nikulin:forms}, the isomorphism classes
of even finite index extensions of a non-degenerate even lattice $M$ are in a canonical bijection 
with the isotropic subgroups ${\Cal G} \subset \discr M$. 
It follows (\cf. \cite[Proposition 1.5.1]{Nikulin:forms}) that
the classes of extensions as in \eqref{eq.primitive-extension} 
are in a bijection with the double cosets 
\[\label{double.cosets}
\gI(M_1) \backslash \{\psi \: \CK \into \discr M_2\} / \gI(M_2),
\]
where $\CK \subset \discr M_1$ is a subgroup and $\psi$ is an injective anti-isometry,
so that the isotropic subgroup ${\Cal G} \subset \discr M_1 \oplus \discr M_2$ as above is the graph of $\psi$.

In the special case where $M$ is a lattice with an involution and $M_1$, $M_2$
are the eigenlattices, $\CK$ is a group of exponent $2$. 
Conversely, if $\CK$ is of exponent $2$, then ${\id_1} \oplus {-\id_2}$
extends to an involution of $M$,
see \cite[Corollary 1.5.2]{Nikulin:forms}.

Back to the proof of \autoref{th.main}, the construction of the prospective equivariant transcendental lattices $T$ from  
a given pair $(T_+, T_-)$
reduces to the computation of the subgroups $\gI(T_\pm) \subset \Aut(\discr T_\pm)$. 
At the same time, for step \iref{item.cosets}, we find the image $\gI_\ttheta(T)$ of the subgroup $O_\ttheta(T) \subset O(T)$
centralizing $\ttheta_T$: it is the subgroup of 
$$
\gI(T_+) \times \gI(T_-) \subset \Aut(\discr T_+) \times \Aut(\discr T_-)
$$
preserving $\CK$ as a set and commuting with $\psi$.

\remark\label{remark:size_K}
In view of \eqref{double.cosets.2} below,
the exponent $2$ subgroup $\CK \subset \discr T_+$ in \eqref{double.cosets}
is subject to an extra restriction
$$
 \ls| \CK |^2 \cdot \ls|\discr\ts_h| = \ls|\discr T_+| \cdot \ls|\discr T_-|.
$$
\endremark

The computation of $\gI(T_\pm)$ is straightforward if the eigenlattice is definite of rank $1$, and is easily done using
Miranda-Morrison theory  \cite{Miranda.Morrison:book} if the eigenlattice is hyperbolic of rank at least $3$.
If it is hyperbolic of rank $2$, we compute $O(T_\pm)$ explicitly.
Up to rescaling, we encounter but the following three classes of lattices
(abbreviating $[a, b, c]:= \Z u + \Z v$, where $u^2 = a, u \cdot v = b, v^2 = c$):
\roster*
\item
$T_\pm \simeq[a,b,0]$
represents~$0$ (equivalently, $-\det T_\pm$ is a perfect square):
the group $O(T_\pm)$ equals $(\Z/2)^2$ or~$\Z/2$
depending on whether the two isotropic
directions can or cannot be interchanged;
\item
$T_\pm \simeq[1,0,c]$: the group $\OG(T_\pm)$ is given by the solutions to
Pell's equation;
\item
$T_\pm \simeq[-2,b,-2]$, $b \geq 3$ odd:
since the generators $u, v$ constitute the two walls of a fundamental polyhedron,
the group $O(T_\pm)$ is generated by~$-1$, the symmetry
$u\leftrightarrow v$ of the polyhedron, and reflection $r_v$,
see~\eqref{eq.reflection}.
\endroster

\step{item.cosets}
Similar to \eqref{double.cosets} in step (\ref{item.cc.T}),
the extensions are in a bijection with the double cosets
\[\label{double.cosets.2}
\gI_{h, \theta}(\ts_h) \backslash \{\psi \: \discr\ts_h \to \discr T\} / \gI_\ttheta(T),
\]
where
\roster*
\item $\psi$ is a \emph{bijective} anti-isometry (since the resulting lattice $\bL$ is unimodular;
any extension \eqref{double.cosets.2}
is isomorphic to $\bL$ as it is unique in its genus),
\item $\gI_\ttheta(T)$ is the subgroup computed in step (\ref{item.cc.T}) together with $T$, and
\item $\gI_{h, \ttheta}(\ts_h)$ is the image of the subgroup $O_{h, \ttheta}(\ts_h) \subset O(\ts_h)$
preserving $h$ and centralizing $\ttheta_S$: it is easily found as a subgroup of the finite group $O(S)$;
in fact, it suffices to use the group of symmetries of the Dynkin diagram.
\endroster

\step{item.sanity.check}
Once a quadruple $(\bL \supset \ts_h \ni h, \ttheta)$
has been constructed,
the verification of condition \iref{i.Nikulin} in \autoref{th.real-statement}
is straightforward.
\qed

\subsection{Proof of \autoref{add.J10}}\label{proof.J10}
Since a sextic curve $C$ in question must have at least two (conjugate)
non-simple singular points, these points are triple, hence adjacent to $\bJ_{10}$,
see \cite{AVG1}. Then, the $\delta$-invariant $\delta(C) \geq 12$,
and $C$ splits into at least three components. Thus, $C$ is a union
of three conics, all passing through a fixed pair of conjugate points $O$, $O'$
with fixed tangent lines $OP$, $O'P$, where $P \in \Rp2$.
It follows that the three conics are members
of a pencil $\Cal P$ with two double base points.
Up to real projective transformation, such a pencil is unique:
it is determined by the triple $(O, O', P)$.
The two degenerate fibers of $\Cal P$, \viz. the double line $OO'$
and the union of $OP$, $O'P$, divides ${\Cal P}_\R$ into two intervals; 
in one of them the conics are empty, whereas in the other their real parts
form a family of nested ovals. This description makes the statement of \autoref{add.J10} immediate.
\qed

\appendix

\section{Symmetric sextics}\label{S.trigonal}

In this appendix, we
describe an explicite geometric construction for the majority of special ($d>1$) empty sextics by means of
the double covering $p\: \Cp2 \dashrightarrow \Sigma_2$
and trigonal curves in the Hirzebruch surface~$\Sigma_2$.
In the extremal cases, where the deck translation of $p$ is a stable involution (see \autoref{ss.sym})
of the sextic in question, this construction is canonical, thus providing also a deformation classification.
As a by-product, we visualize the distinctions between multiple deformation families
within the same real lattice type.

\subsection{Preliminaries}\label{s.trigonal}

Assume that an empty sextic $C\subset\Cp2$
is preserved by an equivariant
involution
$s\:\Cp2\to\Cp2$.
As is known, $s$ has an isolated real fixed point $O\notin C$ (since
$C_\R=\varnothing$)
and a real fixed line
$L$, and the quotient by~$s$ of the blow-up $\Cp2(O)$  is the Hirzebruch
surface $\Sigma_2$ (see \cite{degt:symmetric,degt:8a2,degt:Oka3} for details).
Thus,
we have a commutative diagram
\[
\CD
X&\ {}\overset{\tp}\longdash{}\ &Y\\
@V\pi VV@VV\bar\pi V\\
\Cp2&\ {}\overset{p}\longdash{}\ &\Sigma_2\rlap,
\endCD
\label{eq.diagram}
\]
where each surface is real and each arrow is a (birational)
real double covering. The ramification loci
are as follows:
\roster*
\item
for~$\pi$, the original sextic~$C=p\1(\B)$,
\item
for~$p$, the exceptional section~$E$ and the section $\sec:=p(L)$
disjoint from~$E$,
\item
for~$\barpi$, the exceptional section~$E$ and the image $\B:=p(C)$;
this image is a \emph{proper} (\ie, disjoint from~$E$) trigonal curve in
$\Sigma_2$,
\item
for~$\tp$, the pull-back $\barpi\1(\sec)$.
\endroster
Since the line $L$ has real points, so does~$\sec$ and the real
structure on~$\Sigma_2$ is standard.

To avoid excessive notation, below we systematically use $\bar{\phantom{o}}$ to denote the image under $p$
of a curve in $\Cp2$. Conversely, for a curve ${\bar B} \subset \Sigma_2$, we silently denote by $B$
its pullback in $\Cp2$.

\subsubsection{Trigonal curves and real forms}\label{ss.dessins}
Below, we always start with an extremal trigonal curve $\B\subset\Sigma_2$,
\ie, $\mu(\B)=8$; such a curve appears from a stable involution of a sextic $C$.
(For irreducible sextics this fact is proved
in~\cite{degt:symmetric}; in general, the assertion would follow from
comparing the dimensions of the moduli spaces, but we do not engage into this
discussion and merely state the fact.)
Up to automorphism of~$\Sigma_2$, there are
but finitely many such curves: complex curves are classified by means
of \emph{dessins d'enfants} $\Gamma\subset S^2$,
and the real forms of each curve are the reflections of~$S^2$
preserving~$\Gamma$.

\remark\label{rem.isotrivial}
In \autoref{s.E6+A2}--\autoref{s.E8}, the curve $\B=\B_r$ depends on a
parameter $r\in\R$, so that $\B_0$ is isotrivial whereas all
curves with $r\ne0$ are non-isotrivial and pairwise isomorphic over~$\C$.
In all other cases, the
\emph{complex} curve
with the prescribed set of singularities is unique.
\endremark

\subsubsection{Coordinates}\label{ss.coord}
Till the rest of this appendix, we use real affine coordinates $(x,y)$
in~$\Sigma_2$ so that the exceptional section~$E$ is given by $y=\infty$. In
these coordinates, a proper real section is given by
\[
y=f(x):=ax^2+bx+c,\quad a,b,c\in\R,
\label{eq.section}
\]
and a proper real trigonal curve is given by
\[*
y^3+y^2p_2(x)+yp_4(x)+p_6(x)=0,\quad p_d\in\R[x],\ \deg p_d\le d.
\]
The coordinates about the fiber $x=\infty$ are $u:=1/x$ and $v:=y/x^2$, in
which \eqref{eq.section} takes the form
\[
v=a+bu+cu^2.
\label{eq.infty}
\]

\subsubsection{The construction}\label{ss.construction}
Conversely, given a proper trigonal curve~$\B\subset\Sigma_2$ and a proper
section $\sec\subset\Sigma_2$, the double covering of~$\Sigma_2$ ramified
at $\sec\cup E$ blows down to~$\Cp2$ and the pull-back of~$\B$ is a sextic
curve.
If both~$\B$ and~$\sec$ are real, so is~$C$.

\remark\label{rem.sing}
The set of singularities of~$C$
depends on that of~$\B$ and on the position of~$\sec$ with respect to~$\B$,
see~\cite{degt:symmetric,degt:Oka3}. If
\roster
\item\label{i.generic}
$\sec$ is \emph{generic}, \ie, transverse to~$\B$,
\endroster
then each singular point of~$\B$ doubles in~$C$; thus,
if $\B$ is also extremal, we have $\mu(C)=16$.
Besides, we consider but the following types of degenerate sections:
\roster[\lastitem]
\item\label{i.tangent}
$\sec$ is simple tangent to~$\B$ at a pair of conjugate points, or
\item\label{i.singular}
$\sec$ passes through a pair of conjugate type~$\bA$ singular points of~$\B$.
\endroster
Both degenerations are of codimension~$2$ and result
in sextics~$C$ with $\mu(C)=18$.
\endremark

In the real case, each of the two
complements below admits a chessboard colorings, thus splitting into
two (open) regions:
\[*
(\Sigma_2)_\R\sminus(\B_\R\cup E_\R)=\Sigma_2^+\cup\Sigma_2^-,\qquad
(\Sigma_2)_\R\sminus(\sec_\R\cup E_\R)=\Sigma_2^{+L}\cup\Sigma_2^{-L},
\]
and $C_\R$ is empty if and only if
\[
\text{$\B_\R$ lies in one of the halves~$\Sigma_2^{\pm L}$,
 henceforth denoted by $\Sigma_2^{-L}$},
\label{eq.nonsep}
\]
and the lift of the real structure is chosen so that $\Rp2$ projects
to the closure of $\Sigma_2^{+L}$.
A section satisfying~\eqref{eq.nonsep} is called \emph{non-separating} (with
respect to~$\B$).

Note that \eqref{eq.nonsep} implies (but, in general, is not equivalent to)
that
\[
\text{$\sec_\R$ lies in one of the halves~$\Sigma_2^{\pm}$,
 henceforth denoted by $\Sigma_2^{-}$};
\label{eq.half}
\]
furthermore, the lift of the real structure on~$\Sigma_2$ to~$Y$
in~\eqref{eq.diagram} is to be chosen so that $Y_\R$ project to the closure
of $\Sigma_2^+$.

\lemma\label{lem.section}
For a proper real trigonal curve $\B\subset\Sigma_2$, the space of
non-separating real sections has exactly two connected components
$\CS_\pm(\B)$, which are both convex in the affine space~\eqref{eq.section}
of proper sections.
\endlemma

\proof
The set of non-separating sections is the union of two
open convex sets, \viz.
\[*
\CS_\epsilon:=\bigl\{\text{$f$ as in~\eqref{eq.section}}\bigm|
 \text{$\epsilon f(x)>\epsilon y$ for each point $(x,y)\in\B_\R$}\bigr\},
 \quad\epsilon=\pm.
\]
(At $x=\infty$ this condition is to be modified
according to~\eqref{eq.infty}.)
To prove that, say, $\CS_+\ne\varnothing$, let $b=0$.
By the compactness of $\B_R\sminus E_\R$,
the condition above holds
\roster*
\item
for all $\ls|x|\le1$ whenever $a\ge0$ and $c\gg0$, and
\item
for all $\ls|x|\ge1$ whenever $a\gg0$ and $c\ge0$, \cf.~\eqref{eq.infty}.
\endroster
Thus, $\CS_+$ contains a section
$y=f(x)$ as in~\eqref{eq.section}
with $b=0$ and $a,c\gg0$.
\endproof

\subsubsection{The deformation classification}\label{ss.deform}
A real trigonal curve $\B\subset\Sigma_2$ is said to be
\emph{symmetric} if it is equisingular deformation equivalent,
over~$\R$, to its image under
at least one of the automorphisms
\[*
r\:(x,y)\mapsto(-x,y)\quad\text{or}\quad(x,y)\mapsto(-x,-y)
\]
reversing the orientation of the real fibers of the ruling.
If this is the case, then $r$, followed by the deformation,
interchanges $\Sigma_2^\pm$, as well as
$\CS_\pm(\B)$ in \autoref{lem.section}.

\corollary[of \autoref{lem.section}]\label{cor.deform}
Let $\B\subset\Sigma_2$ be an extremal real trigonal curve, and let $\CM(\B)$
be the space of pairs $(C,s)$, where $C\subset\Cp2$ is an empty sextic,
$\mu(C)=16$, and $s\:\Cp2\to\Cp2$ a real stable involution such that
$C/s\cong\B$. Then $\CM(\B)$ has one \rom(if $\B$ is symmetric\rom)
or two \rom(otherwise\rom) connected components.
\endcorollary

\proof
As explained in \autoref{ss.construction}, the condition $\mu(C)=16$ is
equivalent to the assertion that the section~$\sec$ in the construction
of~$C$ is generic. By \autoref{lem.section}, the space of all sections
satisfying~\eqref{eq.nonsep} has two convex components, which are interchanged by an
automorphism if~$\Sigma_2$ if and only if $\B$ is symmetric, and it remains
to observe that non-generic sections satisfying~\eqref{eq.nonsep} constitute
the intersection with~$\CS(\B)_\pm$ of a real algebraic variety of codimension
at least~$2$, \cf. \autoref{rem.sing}\iref{i.tangent} and~\iref{i.singular}.
\endproof

\subsubsection{Splitting conics and sections}\label{ss.splitting}

Fix a proper trigonal curve
$\B\subset\Sigma_2$. A \emph{splitting section} is a proper
section $\barU\subset\Sigma_2$ such that
\roster
\item\label{i.even}
the local intersection index of $\barU$ and~$\B$ is even at each point, and
\item\label{i.disjoint}
the proper transforms of $\barU$ and~$\B$ in the minimal resolution of
singularities of~$\B$ are disjoint.
\endroster
A proper section $\barU$ is splitting if and only if the proper preimage
$\barpi\1(\barU)$ splits into a pair of disjoint sections of the rational Jacobian
elliptic surface~$Y$, see~\eqref{eq.diagram}.
If $\sec$ is another proper section,
so that we have diagram~\eqref{eq.diagram},
the section $\barU$ is splitting if and only if $B:=p\1(\barU)$ is an
$s$-invariant splitting conic for~$C$.

\remark\label{rem.splitting.section}
Strictly speaking, the existence of an involutive stable symmetry~$s$ and the
presence of splitting conics are two independent properties of the complex lattice 
type; however, often (although not always), the former implies the latter.
Typically, an $s$-invariant splitting conic~$B$ projects to a splitting
section $\barU$. We refrain from a general
statement and merely indicate splitting sections in the models, referring
implicitly to the classification of complex lattice types to identify the
value of~$d$.
\endremark

\observation\label{obs.splitting}
Consider an $s$-invariant real splitting conic~$B$ and
its image $\barU := p(B)$. Due to condition~\iref{i.even} above,
$\barU_\R$ lies entirely in (the closure of) one of the two
regions~$\Sigma_2^\pm$.
It is immediate from the construction that
\roster*
\item
if $\barU$ is in the closure of~$\Sigma_2^-$, see~\eqref{eq.half}, the two
components of $\barpi\1(\barU)$ are complex conjugate and, hence,
$B$ contributes to $c\counts$ in $\scount$ (see \autoref{s.real.invariants});
\item
if $\barU$ is in the closure of~$\Sigma_2^+$, both
components of $\barpi\1(\barU)$ are real and, hence,
$B$ contributes to $r\counts$ in $\scount$.
\endroster
This observation is used to identify the real lattice types obtained in
\autoref{ss.construction}.
\endobservation

\subsubsection{Perturbations of trigonal curves}\label{ss.perturb}
The description of non-isotrivial trigonal curves by means of
\emph{dessins d'enfants} implies that, just like in the case of plane
sextics, {\em all $\bA$-type singular points of a non-isotrivial trigonal curve
$\B\subset\Sigma_2$ can be perturbed arbitrarily and independently.}
We use this observation to realize some real lattice types not admitting an
involutive stable symmetry.

\subsection{The computation}\label{s.computation}
In the rest of this appendix, we consider, one-by-one, the extremal trigonal
curves $\B\subset\Sigma_2$ appearing from the stable involutions~$s$ of the
complex lattice types listed in Tables~\ref{tab.18} and~\ref{tab.16} and, for
each such curve (designated by its sets of singularities),
list its real forms, \cf. \autoref{ss.dessins}.

For each real form, we compute the connected components of the 
equisingular equivariant moduli spaces of pairs $(C,s)$, where $C$ is an empty sextic 
obtained by the construction of \autoref{ss.construction} from~$\B$ and a
section~$\sec$ as in \autoref{rem.sing}. If $\sec$ is generic, we refer to
\autoref{cor.deform}; otherwise, we describe the $1$-parameter family of
sections explicitly, still arriving at two connected components (interchanged
by an automorphism and deformation if $\B$ is symmetric, see \autoref{ss.deform}).

\remark\label{rem.distinguished}
Usually, a real lattice type admits at most one stable involution and,
hence, the moduli space of pairs $(C,s)$ is that of sextics~$C$. In the two
exceptional cases (see Remarks~\ref{rem.2A3+2A1},~\ref{rem.2D4} below), we
explain that one of the three involutions is distinguished and it is this
involution that is used in the construction.
\endremark

By computing the invariants (most notably,
$\ccount$ and $\scount$,
see
\autoref{s.real.invariants}), we establish that each real form of
the given complex lattice type has indeed been realized; 
in most cases, we reestablish the deformation classification 
by showing that each connected component belongs to its own
real form. (Usually, it suffices to compare the counts.)
The
two exceptional cases are emphasized in
\autoref{s.A7+A1} and \autoref{s.D8}
below.

Occasionally, we consider a few perturbations of the original extremal
trigonal curve
$\B\subset\Sigma_2$ (see \autoref{ss.perturb}) and use them to construct
representatives of some other real lattice types, not admitting a stable
involution.

\remark\label{rem.non-extremal}
It is worth
pointing out
that, since, upon the perturbation, the trigonal curve is no longer extremal,
the involution is not stable and
the construction of \autoref{ss.construction} does \emph{not} give us the
complete stratum of sextics. 
\endremark


\subsection{The
trigonal curve $\B(\singset{4A2})$}\label{s.4A2}
The curve~$\B$ is given by
\[
4y^3-(24x^3+3)y+(8x^6+20x^3-1)=0,
\label{eq.4A2}
\]
see~\cite{degt:8a2}; its four cusps are
\[*
P_1=\bigl(0,-\tfrac12\bigr),\quad
 P_{2,3,4}=\bigl(\epsilon,\tfrac32\bigr),\quad \epsilon^3=1.
\]
The non-separating real sections passing through~$P_3$, $P_4$ are
\[
a(x^2 + x + 1) + \frac32,\quad\ls|a+1|>1;
\label{eq.2A5+4A2}
\]
they give rise to the set of
singularities $\sset{2A5+4A2}$.
A generic section gives rise to $\sset{8A2}$.
All sextics are of torus type ($d=3$); this fact follows from the
presence of the splitting sections $\barU_i$ passing through $P_j$, $j\ne i$,
see~\cite{degt:8a2} and \autoref{ss.splitting}.

In spite of its appearance, the curve~$\B$ is symmetric: in appropriate
coordinates, with the cusps at $\pm(1,1)$ and $\pm(2i-i\sqrt3,0)$, it can be
given by a polynomial that is skew-invariant under $(x,y)\mapsto (-x, -y)$.
Hence, in
each of the two cases, we obtain a single real lattice type and a
single deformation family, see \autoref{ss.deform}.

\subsubsection{Perturbations}\label{ss.4A2}
Perturbing the cusp~$P_1$ to $\bA_1$ or $\varnothing$ (see \autoref{ss.perturb}),
we realize the sets of singularities
\[*
\sset{2A5+2A2+2A1},\quad\sset{2A5+2A2},\quad\sset{6A2+2A1},\quad\sset{6A2},
\]
all with $d=3$, as the real splitting section $\barU_1$ remains intact.
The resulting curve is no longer symmetric; hence,
in each of the four cases, we obtain two real forms,
which differ by $\scount$, see
\autoref{s.real.invariants}\iref{(0,1,0)-(1,0,0)} and
\autoref{obs.splitting}.

\subsection{The
trigonal curve $\B(\singset{2A3+2A1})$}\label{s.2A3+2A1}
The curve~$\B$ splits into the zero section $\BD_0\=y=0$ and two other
sections $\BD_\pm$ tangent to~$L_0$ at points $P_\pm$, respectively, and
intersecting at two other points $Q_1,Q_2$. The three real forms  are
(listing $\BD_\pm$ only)
\begin{alignat}4
\label{eq.2A3+2A1.1}
&y=(x\pm1)^2:\quad&&P_\pm=(\pm1,0),\quad&Q_1&=(0,1),\ &Q_2=(\infty,1),\\
\label{eq.2A3+2A1.2}
&y=\pm(x\pm1)^2:\quad&&P_\pm=(\pm1,0),\quad&Q_i&=(\pm i,\pm2i),\\
\label{eq.2A3+2A1.3}
&y=(x\pm i)^2:\quad&&P_\pm=(\pm i,0),\quad&Q_1&=(0,-1),\ &Q_2=(\infty,1).
\end{alignat}
The forms \eqref{eq.2A3+2A1.2} and \eqref{eq.2A3+2A1.3} are symmetric, whereas
\eqref{eq.2A3+2A1.1} is not.
Non-separating real sections passing through $Q_1,Q_2$ exist
in~\eqref{eq.2A3+2A1.2} only; they are
\[
\sec\:\ y=a(x^2+1)+2x,\quad \ls|a|>1,
\label{eq.6A3}
\]
giving rise to the set of singularities
$\sset{6A3}$.
Non-separating real sections
through both $P_\pm$ are
\[
\sec\:\ y=a(x^2+1),\quad\ls|a|>1,
\label{eq.2A7+4A1}
\]
in~\eqref{eq.2A3+2A1.3}; they give rise to the set of
singularities $\sset{2A7+4A1}$.
Finally,
if $\sec$ is a generic non-separating section, we
obtain
$\sset{4A3+4A1}$ (see \autoref{ss.4A3+4A1} below).

In all cases, we have $d=4$ due to the presence of a pair of splitting
sections (not necessarily real) $\barU_i$
passing through both $P_\pm$ and $Q_i$, $i=1,2$.

\remark\label{rem.2A3+2A1}
The sextic~$C$ with the set of singularities $\sset{2A7+4A1}$ has three
stable involutions, and we choose the only one fixing the two type~$\bA_7$
points.
\endremark

\subsubsection{The
sextic $C(\singset{4A3+4A1})$ with $d=4$
and perturbations thereof}\label{ss.4A3+4A1}
The four real forms of
this complex lattice type can be distinguished as follows,
\cf. \autoref{s.real.invariants}\iref{4.forms}.

If $\B$ is as in~\eqref{eq.2A3+2A1.1}, then
$\ccount=(3,0)$.
Both splitting sections
\[*
\barU_i\:\ y=\pm(x^2-1),\quad i=1,2,
\]
are real (hence so are the splitting conics~$B_i$), and their real parts
are in the same half $\Sigma_2^\pm$. Thus, $\scount=(0,2,0)$ or $(2,0,0)$
(see \autoref{obs.splitting}).

If $\B$ is as in~\eqref{eq.2A3+2A1.2}, then
$\barU_1$, $\barU_2$ are conjugate and $\ccount=(3,0)$, $\scount=(0,0,1)$.


If $\B$ is as in~\eqref{eq.2A3+2A1.3}, then
$\ccount=(1,1)$.
Both splitting sections
\[*
\barU_i\:\ y=\pm(x^2+1),\quad i=1,2,
\]
are real and their real parts are in the opposite halves~$\Sigma_2^\pm$;
hence, $\scount=(1,1,0)$.

In the first case,~\eqref{eq.2A3+2A1.1}, we can perturb $Q_2$ (see
\autoref{ss.perturb}),
leaving a single splitting section $\barU_1$, hence still $d=4$.
This gives rise to two real forms,
distinguished by~$\scount$
(see \autoref{s.real.invariants}\iref{(0,1,0)-(1,0,0)} and \autoref{obs.splitting}),
of the set of singularities
$\sset{4A3+2A1}$. The same pair of real forms is obtained by
perturbing~$Q_2$ in the last case~\eqref{eq.2A3+2A1.3}.

In the second case,~\eqref{eq.2A3+2A1.2}, we can perturb a conjugate pair of
nodes of the \emph{original sextic~$C$}.
(This perturbation is no longer symmetric.)
This operation destroys both splitting conics, resulting in the set of
singularities $\sset{4A3+2A1}$ with $d=2$. This real lattice
type is also obtained
from \sset{2D6+2A3}
(see the end of \autoref{ss.D6+2A1.1} below), 
where
$d=2$ in the first place.

\subsubsection{Other perturbations}\label{ss.2A3+2A1}
Perturbing $P_2$ to $\bA_2$ or $\varnothing$
(see \autoref{ss.perturb})
in~\eqref{eq.2A3+2A1.2} and
taking for~$\sec$ either (an appropriate perturbation of)~\eqref{eq.6A3}
or a generic section, we obtain
reducible ($d=2$) curves with the sets of singularities
\[*
\sset{4A3+2A2},\quad \sset{2A3+2A2+4A1},\quad \sset{4A3}.
\]

Perturbing $Q_1,Q_2$ in~\eqref{eq.2A3+2A1.3} and taking
a perturbation of~\eqref{eq.2A7+4A1}
for~$\sec$, we obtain a reducible
($d=2$)
curve
with the set of singularities
$\sset{2A7}$.


\allowdisplaybreaks

\subsection{The
trigonal curve $\B(\singset{2A4})$}\label{s.2A4}
The curve~$\B$ is given by
\begin{gather*}
4y^3-3yp(x)+q(x)=0,\ \text{where}\\
p(x)=x^4-12x^3+14x^2+12x+1,\\
q(x)=(x^2+1)(x^4-18x^3+74x^2+18x+1).
\end{gather*}
see~\cite{degt:Oka3}; it has $\bA_4$ singular points at $(0,1/2)$ and
$(\infty,1/2)$.
A generic non-separating section~$\sec$ gives rise to
the two
real forms for
the set of singularities $\sset{4A4}$.
The curve is special ($d=5$) due to the presence of the splitting sections
\[*
\barU_\pm\:\ y=\frac12(x^2+1)\pm3x.
\]
Both $\barU_\pm$ are real and their real parts are in the same half
$\Sigma_2^\pm$; hence, $\scount=(0,2,0)$ or $(2,0,0)$, see
\autoref{s.real.invariants}\iref{(0,2,0)-(2,0,0)} and
\autoref{obs.splitting}.

\subsection{The
trigonal curve $\B(\singset{A5+A2+A1})$}\label{s.A5+A2+A1}
The curve~$\B$ splits into a ``parabola'' $\BD_2$ and
a section~$\BD_1$ inflection
tangent to~$\BD_2$ at $(1/4,1/2)$:
\[*
\BD_2:\ y^2=x,\qquad \BD_1:\ y=-x^2+\frac32x+\frac3{16},
\label{eq.A5+A2+A1}
\]
see~\cite{degt:2a8}. Taking for~$\sec$ a generic non-separating section, we
arrive at the two real forms for
the set of singularities $\sset{2A5+2A2+2A1}$ with
$d=6$. Indeed, the curve is clearly reducible, hence $2\mathrel|d$, and it is of
torus type
due to the splitting section
\[*
\barU\:\ y=x+\frac14,
\]
resulting also in $\scount=(0,1,0)$ or $(1,0,0)$, see
\autoref{s.real.invariants}\iref{(0,1,0)-(1,0,0)} and
\autoref{obs.splitting}.


\subsubsection{Perturbations}\label{ss.A5+A2+A1}
Perturbing the cusp at infinity
to $\bA_1$ or $\varnothing$ (see \autoref{ss.perturb})
and thus
destroying the torus structure, we arrive at reducible ($d=2$) curves with
the sets of singularities
$\sset{2A5+4A1}$ and $\sset{2A5+2A1}$.

\subsection{The
trigonal curve $\B(\singset{A7+A1})$}\label{s.A7+A1}
The curve~$\B$ splits into a ``conic'' $\BD_2$ with a node at $(\infty,0)$ and a
section $\BD_1$ quadruple tangent to~$\BD_2$ at $(0,1)$:
\[*
\BD_2\:\ \epsilon x^2+y^2=1,\qquad
\BD_1\:\ y=-\frac12\epsilon x^2+1,\quad
 \epsilon=\pm1.
\]
Since both real forms are asymmetric,
a generic non-separating section~$\sec$ gives rise to \emph{two real families}
for
each of the two real forms of
$\sset{2A7+2A1}$.
One has $d=4$
and $\scount=(0,1,0)$ or $(1,0,0)$, see
\autoref{s.real.invariants}\iref{(0,1,0)-(1,0,0)},
due to the splitting section
$\BU\=y=1$.


\subsection{The
trigonal curve $\B(\singset{A8})$}\label{s.A8}
The curve~$\B$ is given by
\[
-y^3+y^2-x^3(2y-x^3)=0
\label{eq.A8}
\]
or, parametrically, by
\[*
x=\frac{t}{t^3+1},\qquad y=\frac1{(t^3+1)^2},
\]
see~\cite{degt:2a8}. According to \cite[Lemma 2.3.3]{degt:2a8}, a real
section~$\sec$ tangent to~$\B$ at two complex conjugate points is uniquely
determined by the tangency points, which must be at
\[
t=-2^{-4/3}\pm si,\quad s>0,\quad 8s^2\ne3\cdot2^{1/3},
\label{eq.2A8+2A1}
\]
and the explicit equation found in~\cite{degt:2a8} shows that $\sec$ is
non-separating, thus giving rise to the set of singularities
$\sset{2A8+2A1}$. A generic section gives rise to $\sset{2A8}$.
In both cases, the curve is of torus type ($d=3$) due to
the splitting section
$\BU\=y=0$.
As $\B$ is asymmetric,
each complex lattice type has two real forms: they differ by
$\scount=(0,1,0)$ or $(1,0,0)$, see
\autoref{s.real.invariants}\iref{(0,1,0)-(1,0,0)} and \autoref{obs.splitting}.

\subsection{The
trigonal curve $\B(\singset{2D4})$}\label{s.2D4}
The curve~$\B$ is isotrivial, splitting into three sections $\BD_0\=y=0$
and $\BD_\pm$. We consider the two real forms
\[*
\BD_\pm\:\ y=\pm x\quad\text{or}\quad y=\pm ix.
\]
A generic section~$\sec$ gives rise to the set of singularities $\sset{4D4}$,
with the two real forms distinguished by
$\ccount=(1,1)$ or $(3,0)$, see \autoref{s.real.invariants}\iref{2.forms}.

\remark\label{rem.2D4}
The curve~$C$ has three stable involutions, and we chose the only one whose action
on $\Sing C$ coincides with that of~$\cc$. Then, we
can assume that all singular points of~$\B$ are real and thus avoid
considering its other real forms, \viz. those with
a pair of complex conjugate singular
points.
\endremark

\subsection{The
trigonal curve $\B(\singset{D5+A3})$}\label{s.D5+A3}
The curve~$\B$ splits into a ``parabola'' $\BD_2$ and section~$\BD_1$
tangent to~$\BD_2$ at $(1,1)$ and passing through its cusp $(\infty, 0)$:
\[*
\BD_2:\ y^2=x,\qquad \BD_1:\ y=\frac12(x+1);
\label{eq.A5+A2+A1}
\]
it has a $\singset{D5}$ singularity at $(\infty,0)$.
Taking for~$\sec$ a generic non-separating section, we
arrive at the two
real forms for the set of singularities $\sset{2D5+2A3}$.
One has $d=4$
and $\scount=(0,1,0)$ or $(1,0,0)$, see
\autoref{s.real.invariants}\iref{(0,1,0)-(1,0,0)},
due to the splitting section $\BU\=y=1$.



\subsection{The
trigonal curve $\B(\singset{D6+2A1})$}\label{s.D6+2A1}
The curve~$\B$ splits into the section $\BD_0\=y=x$ and
two other sections~$\BD_\pm$ intersecting $\BD_0$ at the point $Q=(\infty,0)$
of their tangency and at two
other points $P_\pm$.
Both real forms are symmetric:
\begin{alignat}2
\label{eq.D6+2A1.1}
\BD_\pm\=y&=\pm1:\qquad&P_\pm&=(\pm1,\pm1),\\
\label{eq.D6+2A1.2}
\BD_\pm\=y&=\pm i:&P_\pm&=(\pm i,\pm i).
\end{alignat}
Only \eqref{eq.D6+2A1.2} admits non-separating sections passing
through~$P_\pm$; they are
\[
\sec\:\ y=a(x^2+1)+x,\quad 2\ls|a|>1,
\label{eq.2D6+2A3}
\]
and the resulting sextics have the set of singularities $\sset{2D6+2A3}$.
A generic
section $\sec$ gives rise to the set of singularities $\sset{2D6+4A1}$
(\sstr{\bD_6}).
In both cases, $d=2$.
Hence, $d\divides|2$ (and thus $d=2$, as we keep the curves reducible) for
all perturbations below.

\subsubsection{Perturbations of $\sec$ as in \eqref{eq.2D6+2A3} and
 $\B$ as in~\eqref{eq.D6+2A1.2}}\label{ss.D6+2A1.1}
We can
perturb $Q$ to
\[*
\singset{D4+A1}\to\singset{D4}\quad\text{or}\quad
\singset{A3+2A1}\to\singset{4A1}\overset*\to\singset{3A1}\overset*\to\singset{2A1},
\]
arriving at the sets of singularities
\[*
\aligned
&\sset{2D4+2A3+2A1}\mbox{\ (\str{\bA_3})},\quad\sset{2D4+2A3},\quad\text{or}\\
&\sset{4A3+4A1},\quad
 \sset{2A3+8A1}\mbox{\ (\str{\bA_3}, \cf. \autoref{ss.D6+2A1.2})},\quad
 \sset{2A3+6A1},\quad \sset{2A3+4A1},
\endaligned
\]
respectively, all with $d=2$. Here, the two codimension one perturbations are
easily described explicitly,
\[*
\alignedat2
\BD_\pm\:\ &y=\pm(1-\Ge)+\Ge x\quad&\text{or}\quad &y=\pm i(1-\Ge)+\Ge x,
\rlap{\quad or}\\
\BD_0\:\ &y=x+\Ge(x^2-1)\quad&\text{or}\quad &y=x+\Ge(x^2+1),
\endalignedat
\]
respectively (we describe both real forms of~$\B$), upon which we can use
\autoref{ss.perturb} to perturb type~$\bA$ singular points. In the last two
cases, marked with a~$^*$, we perturb one or both points of intersection
of~$\BD_+$ and~$\BD_-$, thus keeping $D_0$ a separate conic component
of~$C$.

Alternatively, starting from $\singset{A3+2A1}$ and dissolving one of the nodes,
we arrive at the set of singularities $\sset{4A3+2A1}$ with $d = 2$,
\cf. the end of \autoref{ss.4A3+4A1}. 


\subsubsection{Perturbations of $\sec$ generic and
 $\B$ as in~\eqref{eq.D6+2A1.1} or~\eqref{eq.D6+2A1.2}}\label{ss.D6+2A1.2}
The perturbations
\[*
\singset{D4+A1},\quad\singset{A3+2A1},\quad\singset{4A1}
\]
in \autoref{ss.D6+2A1.1} produce two real forms,
with $\ccount=(1,1)$ or $(3,0)$, see
\autoref{s.real.invariants}\iref{2.forms},
for each of
\[*
\sset{2D4+6A1},\quad
\sset{2A3+8A1}\mbox{\ (\sstr{\bA_3}, \cf. \autoref{ss.D6+2A1.1})},\quad
\sset{12A1}.
\]
The last two perturbations result in the sets of singularities $\sset{10A1}$ and $\sset{8A1}$,
both sextics splitting into a quartic and a conic.


\subsubsection{Perturbations of $\sec$ generic and
 $\B$ as in~\eqref{eq.D6+2A1.1}}\label{ss.D6+2A1.3}
In addition to~\autoref{ss.D6+2A1.2}, we can also perturb~$Q$ to
$\singset{D5}$ (perturbing $\BD_+\cup\BD_-$ to a ``parabola''~$\BD$)
or $\singset{A2+2A1}$:
\[*
\BD\:\ (y^2-1)+\Ge x=0,\qquad
\BD_0\:\ y=x+\Gd x^2,
\]
arriving at the sets of singularities
\[*
\sset{2D5+4A1},\quad \sset{2A2+8A1}.
\]
Alternatively, the perturbation of the node $\BD_+\cap\BD_0$, possibly preceded
by
the perturbation of $Q$ to $\singset{D4+A1}$ as in \autoref{ss.D6+2A1.1},
gives rise
to
\[*
\sset{2D6+2A1},\quad\sset{2D4+4A1}.
\]
All four sextics obtained split into a quartic and a conic.

\subsection{The
trigonal curve $\B(\singset{D8})$}\label{s.D8}
The trigonal curve~$\B$ splits into the section $\BD_1=\{y=x\}$ and a
``conic'' $\BD_2$ with a node at $(\infty,0)$; the two real forms are
\begin{align}
\label{eq.D8.1}
\BD_2\=y^2-x^2&=1,\\
\label{eq.D8.2}
\BD_2\=x^2-y^2&=1.
\end{align}
Arguing as in~\cite{degt:8a2,degt:Oka3}, we
conclude that,
if a section $\sec$ is bitangent to $\BD_2$, the two tangency points
must be of the form
$(\pm x_0, y_0)$.
These points can be complex conjugate only in case~\eqref{eq.D8.1},
and a double tangent section $\sec$ does exist if and only if $0<\ls|y_0|<1$.
Then
$\sec$ is automatically non-separating and gives rise to the set of singularities
$\sset{2D8+2A1}$. A generic section
in \eqref{eq.D8.1} or \eqref{eq.D8.2} gives rise to the
set of singularities $\sset{2D8}$.
As $\B$ is symmetric, each real form of~$\B$ results in a single deformation
family; thus, we obtain \emph{two deformation families} for \sset{2D8}.
%
All sextics obtained are reducible, \ie, we have $d=2$.

\subsubsection{Perturbations}\label{ss.D8}
Perturbing the type $\bD_8$ point at infinity
to $\bD_6\oplus\bA_1$ (\eg, replacing~$\BD_1$ with
the line
$y=x+\Ge$) and keeping~$\sec$ bitangent to~$D$,
we obtain the set of singularities \sset{2D6+4A1} (\str{\bD_6},
\cf. \autoref{s.D6+2A1}).

\subsection{The
trigonal curve $\B(\singset{E6+A2})$}\label{s.E6+A2}
The curve~$\B$ is given by
\[*
y^3+(ry+x)^2=0,\quad r\in\R,
\]
see \cite{degt:e6}; it has an $\singset{E6}$ singularity at $(\infty, 0)$,
and a generic non-separating section~$\sec$  gives rise
to the set of singularities $\sset{2E6+2A2}$.
The sextics are of torus type ($d=3$) due to the splitting section
$\BU\=y=0$.
Since $\B$ is
asymmetric, there are two
real forms,
with $\scount=(0,1,0)$ or $(1,0,0)$, see
\autoref{s.real.invariants}\iref{(0,1,0)-(1,0,0)} and
\autoref{obs.splitting}).

\subsection{The
trigonal curve $\B(\singset{E7+A1})$}\label{s.E7+A1}
The curve~$\B$ splits into a ``parabola'' $\BD_2$ and section~$\BD_1$:
\[*
\BD_2:\ y^2=x,\qquad \BD_1:\ y=r,\quad r\in\R;
\label{eq.A5+A2+A1}
\]
it has an $\singset{E7}$ singularity at $(\infty,0)$.
Taking for~$\sec$ a generic non-separating section, we
arrive at the set of singularities $\sset{2E7+2A1}$ with
$d=2$ (reducible curve).
All curves $\B_r$ are symmetric;
hence, we obtain a single
connected deformation family.

\subsection{The
trigonal curve $\B(\singset{E8})$}\label{s.E8}
The curve~$\B$ is parametrized by
\[*
x=t^3+3rt, \quad y=t, 
\]
where $r \in \R$. 
The sign of $r\ne0$ selects one of the two distinct real forms of
non-isotrivial curves, see \autoref{rem.isotrivial}.
Arguing as in~\cite{degt:8a2,degt:Oka3}, we
conclude that a section double tangent to~$\B$
at two values $t_1\ne t_2$ of the parameter
exists if and only if
\[*
r=t_1^2+3t_1t_2+t_2^2.
\]
The values $t_{1,2}$ can be complex conjugate,  $t_{1,2} =\Ga\pm\Gb i$, if and only if
$$r=5\alpha^2+\beta^2>0,$$
and a double tangent section does exist
if and only if $\Ga\ne0$. In this case, the section is automatically
non-separating, giving rise to the set of singularities $\sset{2E8+2A1}$.
A generic section gives rise to $\sset{2E8}$.
Since all curves~$\B_r$ are symmetric, we obtain a single connected deformation
family in each of the two cases.

\section{Further examples}\label{S.explicit}

In this appendix, we consider a couple of examples illustrating the computation 
leading to the proof of \autoref{th.main}. 
In \autoref{s.d=1} we give a simple proof of the fact that the necessary condition 
\eqref{i.embeddable} is also sufficient in the case of non-special sextics. 
In \autoref{s.2A9} we show that, as one would expect,
the two real deformation families with the set of singularities $2\bA_9$
correspond to the two distinct complex ones. 
 
\subsection{The non-special curves}\label{s.d=1}
In the case of non-special ($d = 1$) empty sextics,
the sufficiency of condition \eqref{i.embeddable} (\cf. \autoref{rem.double})
can easily be proved directly,
without a reference to the classification.
%

\theorem\label{th.non-special}
A non-special empty sextic with an even set of singularities $2\ccS$
exists if and only if
the root system $\ccS$ admits a primitive embedding
into
$\bE_8 \oplus \bU$,
or, equivalently, the Dynkin diagram of $\ccS$ is an induced subgraph of \eqref{eq.E10}.
\endtheorem

\proof
As explained in \autoref{th.real-statement} and \eqref{L+L-},
the empty involution on $\bL$ can be described as follows (see \cite{Nikulin:forms}): 
the eigenlattices are
\[*
\alignedat2
L_+&= \bE_8(2) \oplus \bU(2), && \quad \text{with a standard basis $e'_1, \ldots, e'_{10}$, \cf. \eqref{eq.E10}}, \\
L_-&= \bE_8(2) \oplus \bU(2) \oplus \bU,  && \quad \text{with a standard basis $e''_1, \ldots, e''_{10}$, $u_1, u_2$},
\endalignedat
\]
and $\bL$ is the extension of $L_+ \oplus L_-$ \via\/ the vectors $e^\pm_i := \frac12(e'_i \pm e''_i)$, $i = 1, \ldots, 10$.
Now, we assume that $\ccS \subset L_+$ is generated by a subset $e'_i$, $i \in I$,
of the basis, and it is immediate that the vectors $e^\pm_i$, $i \in I$,
generate the set of singularities $2\ccS$.
The polarisation is $h := u_1 + u_2$, and the existence of a sextic
is given by \autoref{th.real-statement}.

The necessity of the condition follows from \autoref{cor.reduction}
and obvious observation that, if $S := 2\ccS$ is primitive in $\bL$,
then so are $S_\pm \subset L_\pm$.
\endproof

\subsection{The set of singularities $\singset{2A9}$}\label{s.2A9}
This is the only complex lattice type whose
corresponding equisingular family has two real components of positive
dimension, see \cite{Aysegul:paper}.
We assert that each of the two components has an
empty real representative. We have $\ts_h=2\bA_9\oplus\Z h$, and the
transcendental lattice is
\[*
T:=\Z u\oplus\Z v_1\oplus\Z v_2,\quad u^2=-2,\quad v_{1,2}^2=10;
\]
the two lattices can be glued, in the sense of \eqref{double.cosets.2}, by the anti-isometry~$\psi$ 
\[*
\gathered
\tfrac12u \mapsto [0,0,1],\\
\tfrac12v_1\mapsto[5,0,0],\qquad\tfrac25v_1\mapsto[2,0,0]\rlap,\\
\tfrac12v_2\mapsto[0,5,0],\qquad\tfrac25v_2\mapsto[0,2,0]
\endgathered
\]
of the discriminants; we use the shorthand notation $[a, b, c]$ for the elements
of
$$
\discr(\bA_9 \oplus \bA_9 \oplus \Z h) = (\Z/10) \oplus (\Z/10) \oplus (\Z/2).
$$
The ``obvious'' automorphisms of~$T$, \viz. inverting
one of the generators and the transposition $v_1\leftrightarrow v_2$,
generate a subgroup $G\subset\Aut\discr T$ of index $2$.
On the other hand, a computation using
\cite{Miranda.Morrison:book} shows that the full image of
$\OG(T)$ in $\Aut\discr T$ is of index~$2$; hence, the latter image
equals~$G$. Furthermore, $\psi_*(G)$ equals the image of $\OG_h(\ts_h)$ in
$\Aut\discr\ts_h$, and it is this fact that gives rise to two
homological types/connected components: $\ts_h$ and~$T$ can be glued \via~$\psi$
or \via\ $\psi':=\Gf\circ\psi$, where $\Gf$ is, say, the transposition
$[2,0,0]\leftrightarrow[0,2,0]$ of the two generators of $\discr_5\ts_h$.
A straightforward computation shows that, in both cases, the involution
\[*
-\bigl((\bA_9\leftrightarrow\bA_9)\oplus(v_1\leftrightarrow v_2)\bigr)
\]
extends to one induced by an empty real structure.



%

{
\let\.\DOTaccent
\def\cprime{$'$}
\bibliographystyle{amsplain}
\bibliography{degt}
}

\end{document}